\documentclass[12pt]{article}
\usepackage{amsfonts}
\usepackage{amsmath}
\usepackage{amssymb}
\usepackage{amscd}

\begin{document}
\begin{center}
{\large \bf  ANALOGIES BETWEEN KNOTS AND PRIMES, 3-MANIFOLDS AND NUMBER RINGS}
\end{center}
\begin{center}
{\em Dedicated to the memory of my mother, Chieko Morishita}
\end{center}
\begin{center}
{Masanori MORISHITA}
\end{center}
{\small ABSTRACT.  This is an expository article of our work on
analogies between knot theory and algebraic number theory. We shall discuss foundational analogies between knots and primes,
  3-manifolds and number rings mainly from the group-theoretic point of
view.}
\footnote[0]{2000 {\em Mathematics Subject Classification.} Primary
11R, 57M}
\footnote[0]{This article is the expanded English translation of [56] with up-dated references and }
\footnote[0]{will appear in Sugaku Expositions, AMS.}
\footnote[0]{The original version [55] of this article was written during 2002-2003.}
\\
\\
{\bf Contents}\\
Introduction \\
1. Knots and primes \\
2. Linking numbers and Legendre symbols\\
3. Link groups and pro-$l$ Galois groups with restricted ramification  \\
4. Milnor invariants and multiple power residue symbols \\
5. Homology groups and ideal class groups\\
6. 3-manifold coverings and number field extensions \\
6.1. Hilbert theory\\
6.2. Capitulation problem \\
6.3. Class towers\\
7.  Alexander-Fox, torsion theory and Iwasawa theory \\
8. Moduli and deformations of representations of knot and prime groups\\
\\
{\bf Introduction}
\\

In this article, we try to bridge two fields, knot theory and algebraic number theory, which branched out after the works of C.F. Gauss and have grown up in their separated ways during a century and a half, from the viewpoint of the analogies between knots and primes. We hope that our attempt would provide new insights and problems to both fields.\\

 Historically, the 3-dimensional view for a number ring was recognized
for the first time when J. Tate, M. Artin and J.-L. Verdier gave a topological interpretation of classfield theory. Namely, the classical classfield theory by T. Takagi - E. Artin is restated as a sort of 3-dimensional Poincar\'{e} duality in the \'{e}tale cohomology of a number ring ([44]). The analogy between knots and primes was firstly pointed out by B. Mazur from this homotopical viewpoint, during that time, in the middle of 1960's (and also by Y. Manin). However, it would be fair to say that even during this period and afterward, number fields have been commonly regarded as analogues of algebraic function fields of one variable along with the development of arithmetical algebraic geometry. In recent years, M. Kapranov and A. Reznikov took up the analogy between  number fields and 3-manifolds again ([29],[70],[71]) and christened the study of those analogies arithmetic topology. My own independent study of the analogies between knots and primes was based on an analogy I noticed between the structures of a Galois group with restricted ramification and of a link group. In particular, the analogies between classical invariants in knot theory and algebraic number theory such as linking numbers and power residue symbols, Alexander polynomials and Iwasawa polynomials are clearly explained from this group-theoretic point of view  ([51],[54],[57],[58],[59]). For instance, we can introduce the arithmetic analogues of Milnor's higher linking numbers and then Gauss' genus theory is generalized in a natural way from the standpoint of this analogy. In this article, I will explain the basic ideas, motivations and results in my recent work, hoping that this survey would provide a foundation for the future study of arithmetic topology.\\

As we shall see in the following, the analogy between primes and knots is well-suited with the arithmetic (for example, classfield theory) originated from Gauss' Disquisitiones Arithmeticae. Furthermore, as Gauss' linking number was connected to the electro-magnetic theory from its origin, the analogy between number theory and knot theory seems to be extended to the connection involving quantum field theory, dynamics etc (cf. [8],[28],[42] etc). I would like to e-mail Gauss in heaven and ask ``Is this circle of thoughts what you dreamed to explore ?".\\

It is my pleasure to thank J. Hillman, M. Kapranov, T. Kohno, B. Mazur, J. Morava, K. Murasugi and T. Ono for useful, enlightening and encouraging discussions and correspondences. This work is partly supported by the Grants-in-Aid for Scientific Research, Ministry of Education, Culture, Sports, Science and
Technology, Japan.\\

Throughout this article, we write $\pi_{1}(X)$ (omitting a base point) for the fundamental group (resp. \'{e}tale fundamental group) of a topological space (resp. a scheme) $X$. A manifold is  assumed to be oriented and connected.\\
\\
{\bf 1. Knots and primes}\\

A circle $S^{1}$ is homotopically the Eilenberg-MacLane space $K({\bf Z},1)$ and hence its arithmetic analog is $K(\hat{\bf Z},1)$ ($\hat{\bf Z}$ being the pro-finite completion of ${\bf Z}$) \'{e}tale homotopically, namely, a finite field ${\rm Spec}({\bf F}_{q})$:
$$ S^{1}\;\; \longleftrightarrow \;\; {\rm Spec}({\bf F}_{q}). \leqno{(1.1)}$$
A finite cyclic covering $S^{1} = {\bf R}/n{\bf Z} \rightarrow S^{1} = {\bf R}/{\bf Z}$ corresponds to a finite cyclic extension ${\bf F}_{q^{n}}/{\bf F}_{q}$, and the loop $l \in \pi_{1}(S^{1}) = {\rm Gal}({\bf R}/S^{1})$, $l(x) = x+1$ ($x \in {\bf R})$), corresponds to the Frobenius automorphism ${\rm Fr} \in \pi_{1}({\rm Spec}({\bf F}_{q})) = {\rm Gal}(\overline{\bf F}_{q}/{\bf F}_{q})$, ${\rm Fr}(x) = x^{q}$ ($x \in \overline{\bf F}_{q})$:
$$ \;\;\;\;\;\;\;\; \in \pi_{1}(S^{1}) \;\; \longleftrightarrow \; \; {\rm Fr} \in \pi_{1}({\rm Spec}({\bf F}_{q})) .$$
A tubular neighborhood $V = S^{1} \times D^{2}$ ($D^{2}$ is a 2-dimensional disk) of $S^{1}$ is homotopy equivalent to $S^{1}$ and $V \setminus S^{1}$ is homotopy equivalent to the 2-dimensional torus $\partial V$. On the other hand, for a $\frak{p}$-adic integer ring ${\cal O}_{\frak{p}}$ with residue field ${\bf F}_{q}$ and the $\frak{p}$-adic field $k_{\frak{p}}$, ${\rm Spec}({\cal O}_{\frak{p}})$ is homotopy equivalent to ${\rm Spec}({\bf F}_{q})$ and ${\rm Spec}({\cal O}_{\frak{p}}) \setminus {\rm Spec}({\bf F}_{q})$ is  ${\rm Spec}(k_{\frak{p}})$. Thus ${\rm Spec}({\cal O}_{\frak{p}})$ and ${\rm Spec}(k_{\frak{p}})$ are seen as analogues of $V$ and $\partial V$ respectively. In fact, we find the analogies between their fundamental groups as follows. For the natural homomorphism $\pi_{1}(\partial V) \rightarrow \pi_{1}(V) = \pi_1(S^1)$, its kernel is the infinite cyclic group generated by a loop $\alpha = \{ a \} \times \partial D^{2}$ ($a \in S^{1}$), called a meridian,  and the inverse image $\beta = S^{1} \times \{ b \}$ ($b \in \partial D^{2}$) of $l \in \pi_1(S^1)$ is called a longitude. The group $\pi_{1}(\partial V)$ is a free abelian group generated by $\alpha$ and $\beta$ which are subject to the relation $[\alpha,\beta]:=\alpha\beta\alpha^{-1}\beta^{-1}=1$. On the other hand, for the natural homomorphism $\pi_{1}({\rm Spec}(k_{\frak{p}})) \rightarrow \pi_{1}({\rm Spec}({\cal O}_{\frak{p}}))=\pi_{1}({\rm Spec}({\bf F}_{q}))$, its kernel $I_{\frak{p}}$ is called the inertia group and the inverse image $\sigma$ of ${\rm Fr} \in \pi_{1}({\rm Spec}({\bf F}_{q}))$ is also called (an extension of) Frobenius automorphism.  The maximal tame quotient of $I_{\frak{p}}$ is topologically generated by a monodromy $\tau$ and the maximal tame quotient $\pi_{1}^{tame}({\rm Spec}(k_{\frak{p}}))$ of $\pi_{1}({\rm Spec}(k_{\frak{p}}))$ is topologically generated by $\tau$ and $\sigma$ which are subject to the relation $\tau^{q-1}[\tau,\sigma]=1$.

$$ \begin{array}{ccc} V & \longleftrightarrow & {\rm Spec}({\cal O}_{\frak{p}}),\\
\partial V & \longleftrightarrow & {\rm Spec}(k_{\frak{p}}),\\
 \mbox{meridian} \;  \alpha  & \longleftrightarrow &\mbox{monodromy}\; \tau,\\
{\rm longitude} \; \beta & \longleftrightarrow & \mbox{Frobenius auto.} \; \sigma,
\\
\pi_{1}(\partial V) = \langle  \alpha, \beta \; | \; [\alpha, \beta] =
1 \rangle & \longleftrightarrow  & \pi_{1}^{tame}({\rm Spec}(k_{\frak{p}})) = \langle 
\tau, \sigma\; | \; \tau^{q-1}[\tau,\sigma]=1 \rangle.
\end{array} \leqno{(1.2)}
$$
\vspace{1.8cm}

A knot is an embedding of $S^{1}$ into a 3-dimensional manifold $M$. On the other hand, for the ring of integers ${\cal O}_{k}$ of a number field $k$ of finite degree over ${\bf Q}$, it was shown by M. Artin and J.-L. Verdier ([44],[47])
that ${\rm Spec}({\cal O}_{k})$ has cohomological dimension 3 (up
to 2-torsion) and enjoys a sort of 3-dimensional Poincar\'{e} duality in
\'{e}tale cohomology. Thus ${\rm Spec}({\cal O}_{k})$ is seen as an analogue of a 3-manifold:
$$ \mbox{3-manifold} \; M\;\; \longleftrightarrow \;\;  {\rm Spec}({\cal O}_{k}). \leqno{(1.3)}$$
For a prime ideal $\frak{p} (\neq 0)$ of ${\cal O}_{k}$ with residue field ${\bf F}_{\frak{p}} := {\cal O}_{k}/\frak{p}$, the natural map ${\rm Spec}({\bf F}_{\frak{p}}) \hookrightarrow {\rm Spec}({\cal O}_{k})$ is then seen as an analog of a knot in view of (1.1) and (1.3):
$$\mbox{knot} \; S^{1} \hookrightarrow M \; \longleftrightarrow \; \mbox{prime}\; {\rm Spec}({\bf F}_{\frak{p}}) \hookrightarrow
 {\rm Spec}({\cal O}_{k}) .$$
Since $\pi_{1}({\rm Spec}({\bf Z}))
= 1$ (Hermite-Minkowski) in particular, by analogy with the Poincar\'{e} conjecture,
${\rm Spec}({\bf Z})$ added by the infinite prime is seen as an analog
of the standard 3-sphere $S^{3}$ and the prime ideal $(p)$ of ${\bf Z}$ ($p$ being a prime number) looks like a knot in ${\bf R}^{3}$:
$$ S^{1} \hookrightarrow {\bf R}^{3} \cup \{ \infty \} = S^{3} \;
\longleftrightarrow \; {\rm Spec}({\bf F}_{p}) \hookrightarrow {\rm Spec}({\bf Z}) \cup \{ \infty \}.$$
Here we regard $S^3$ as the end-compactification of ${\bf R}^3$ by identifying $\infty$ with the end of ${\bf R}^3$, and see the infinite prime of $\bf Q$ as an analogue of the end. In general, the set of infinite primes of a number field may be seen as an analogue of the end of a non-compact 3-manifold ([8],[68]).\\

For a knot $K$ in a closed 3-manifold $M$, the complement $X_{K}$ of the interior $V_{K}^{o}$ of a tubular neighborhood $V_K$ of $K$ in $M$ is a compact 3-manifold with boundary $\partial V_{K}$ and the fundamental group $G_{K} := \pi_{1}(X_{K}) = \pi_{1}(M\setminus K)$ is called the knot group of $K$. The knot group $G_{K}$ reflects how $K$  is knotted in $M$. In fact, it was shown by by W. Whitten, C. Gordon-J. Luecke that for prime knots $K$ and $L$ in $S^{3}$,
$$ G_{K}  \simeq G_L \;  \Longleftrightarrow \; K \simeq L \; (\mbox{up to orientation}).\leqno{(1.4)}$$
A meridian and a longitude of a knot $K$ means the 2 loops $\alpha$ and $\beta$ on $\partial V_{K} = \partial X_{K}$ as in (1.2) respectively. The group $D_{K} := \pi_{1}(\partial X_{K}) = \langle \alpha, \beta | [\alpha,\beta]=1\rangle$ is called the peripheral group of $K$. For the case $M = S^{3}$, the knot group $G_{K}$ has a Wirtinger presentation ([5]) which implies that $G_{K}$ is generated by conjugates of $I_{K} = \langle \alpha \rangle$. \\
\\
On the other hand, according to (1.2), the $\frak{p}$-adic field ${\rm Spec}(k_{\frak{p}})$ for a prime ideal $\frak{p}$ of ${\cal O}_{k}$ plays a role of the ``boundary" of $X_{\frak{p}} := {\rm Spec}({\cal O}_{k}) \setminus \{ \frak{p} \}$, and the fundamental group $G_{\{ \frak{p} \}} := \pi_{1}(X_{\frak{p}})$ reflects how a prime $\frak{p}$ is knotted in ${\rm Spec}({\cal O}_{k})$. Following the case of a knot, we call $G_{\{ \frak{p}\}}$ the prime group of $\frak{p}$. As an analogy of (1.4), we note the following: for prime numbers $p$ and $q$, we have
$$ G_{\{ (p) \}} \simeq G_{\{(q)\}} \; \Longleftrightarrow \; p = q. \leqno{(1.5)}$$
The analogue of the peripheral group $D_K$ is  the decomposition group $D_{\{ \frak{p}\}}:= \pi_{1}({\rm Spec}(k_{\frak{p}}))$ (or its maximal tame quotient for closer analogy). As an analogy of the Wirtinger presentation, we see that $G_{\{(p) \}}$ is topologically generated by conjugates of the inertia group $I_{\{ (p) \}}$:
$$ \begin{array}{ccc}   \partial V_{K} \subset M \setminus V_{K}^{o} & \longleftrightarrow &  {\rm Spec}(k_{\frak{p}})\subset {\rm Spec}({\cal O}_{k}) \setminus \{ \frak{p} \},\\
\mbox{peripheral group}\; D_{K} \rightarrow G_{K}
 & \longleftrightarrow &  \mbox{decomposition group} \; D_{\{ \frak{p} \}} \rightarrow G_{\{ \frak{p} \}}.
\end{array} \leqno{(1.6)}
$$

For a number field $k$ and a finite set $S$ of prime ideals, the \'{e}tale fundamental group $\pi_{1}({\rm Spec}({\cal O}_{k} \setminus S))$, which is the Galois group $G_{S} = {\rm Gal}(k_{S}/k)$ of the maximal extension $k_{S}$ over $k$ unramified outside $S$ and the set of infinite primes, is seen as an analogue of a link group:
$$ \mbox{link group} \; G_{L} = \pi_{1}(M \setminus S)
\longleftrightarrow \mbox{Galois group} \; G_{S} = \pi_{1}({\rm Spec}({\cal O}_{k} \setminus S)).$$
The Galois group $G_{S}$ is huge in general, and it is unknown whether even $G_{\{(p)\}}$ is finitely generated or not. The author thinks that it is one of the  fundamental problems in number theory to understand\\
\\
\hspace{2cm}{\em how a prime ${\rm Spec}({\bf F}_{p})$ is embedded in ${\rm Spec}({\bf Z})$ }\\
\\
and this is the very theme studied in knot theory.\\

It has been commonly known that ${\rm Spec}({\cal O}_{k})$ is an analogue of an algebraic curve $C$ over a finite field ${\bf F}_{q}$. By the exact sequence $1 \rightarrow \pi_1(C\otimes \overline{\bf F}_{q}) \rightarrow \pi_1(C) \rightarrow \pi_1({\rm Spec}({\bf F}_{q})) \rightarrow 1$, a curve $C$ is seen homotopically as a surface-bundle over a ``circle" ${\rm Spec}({\bf F}_{q})$ with fiber $C\otimes \overline{\bf F}_{q}$. Since a number ring has no constant field, there is no analogy between the structures of $\pi_1({\rm Spec}({\cal O}_{k}))$ and $\pi_1(C)$, and in fact $\pi_1({\rm Spec}({\cal O}_{k}))$ behaves in a quite random manner like 3-manifold fundamental groups $\pi_1(M)$. We also note that the analogy of (1.4) and (1.5) does not hold for two points on an affine line. In Iwasawa theory, for a prime number $p$, a ${\bf Z}_{p}$-extension $k_{\infty}/k$ obtained adjoining all $p^n$-th roots of unity for $n\geq 1$ is regarded as an analog of the constant field extension ([26]). However,   $k_{\infty}/k$ is an extension ramified over $p$, while $C\otimes \overline{\bf F}_{q}$ is an unramified covering of $C$. It is our viewpoint that $k_{\infty}/k$ should be regarded as an analog of the tower of cyclic coverings of 3-manifolds branched along a knot (cf. Section 7). Thus, it seemed to me that it is more natural to see the Galois group $G_{S}$ above as an analogue of a link group rather than a fundamental group of a Riemann surface. This was one of the motivations for my study. Note that the viewpoint to see a prime as a loop has been known in the study of closed geodesics on a Riemannian manifold or of closed orbits of a dynamical system ([7],[80]). However, our analogy is proper in 3-dimension in the respect that we consider the local analogies such as (1.2) and (1.6).\\
\\
{\bf 2. Linking numbers and Legendre symbols}
\\

The origins of knot theory and algebraic number theory may be found in the works by Gauss on linking numbers ([14]) and quadratic residues ([13]). When we re-examine both notions from the viewpoint of the analogies in Section 1, we find the close analogy between them. This is another motivation for my work.\\

One way to detect the linking number is to regard it as a monodromy.
To avoid the ambiguity of sign of the linking number caused by the
orientation, we consider the mod 2 linking number. Suppose that $K \cup L$ is a link of two components in $S^{3}$. Take the unique
unbranched double  cover  $Y_{K} \rightarrow X_{K} := S^{3} \setminus
K$. We then get the linking number $\mbox{lk}(K,L)$ mod 2 as the
covering transformation defined by a longitude $\beta_L$ around $L$:
$$ \begin{array}{ccc}
\pi_{1}(X_{K}) &  \longrightarrow  &  \mbox{Gal}(Y_{K}/X_{K}) =  {\bf
Z}/2{\bf Z}  \\

  [\beta_L]  &  \longmapsto &     \;\;\;\;\;\;\;\;
\mbox{lk}(K,L) \; \mbox{mod} \; 2
\end{array} $$

Another way to get the linking number is to use the intersection
number of $K$ with the Seifert surface $\Sigma$ of $L$, $\partial \Sigma = L$.
 Namely, we get the linking number as the cup product of the cohomology
classes $K$ and $\Sigma$: For $X_L := S^3 \setminus L$, 
$$ \begin{array}{ccc}
H^{2}_{c}(X_{L},{\bf F}_{2}) \times H^{1}(X_{L},{\bf F}_{2})  &
\stackrel{\cup}{\longrightarrow} &  H^{3}_{c}(X_{L},{\bf F}_{2}) =  {\bf
Z}/2{\bf Z}\\
([K], [\Sigma]) & \longmapsto & \;\;\;\;\;  [K] \cup [\Sigma] =
\mbox{lk}(K,L) \, \mbox{ mod} \; 2.\\
\end{array}$$

On the other hand, suppose that $p$ and $q$ are distinct odd
  prime numbers. Since the quadratic extension ${\bf Q}(\sqrt{p})/{\bf Q}$ is ramified only over $p$ if and only if $p \equiv 1$ mod $4$, we assume that $p, q \equiv 1$ mod $4$.  Take the unique  \'{e}tale
double covering $Y_{p}  \rightarrow X_{p} := \mbox{Spec}({\bf Z}) \setminus
\{ (p) \} = {\rm Spec}({\bf Z}[\frac{1}{p}])$, where $Y_{q}$ is the spectrum of the normalization of ${\bf Z}[\frac{1}{p}]$ in ${\bf Q}(\sqrt{p})$.  According to (1.2), we may define the
mod 2 linking number of $p$ and $q$, denoted by $\mbox{lk}_{2}(p,q)$, by the image of the conjugacy class of the Frobenius automorphism over $q$ in $\mbox{Gal}(Y_{p}/X_{p}) =
{\bf Z}/2{\bf Z}$:
$$ \begin{array}{ccc}
 \pi_{1}(X_{p}) &  \longrightarrow &  \mbox{Gal}(Y_{p}/X_{p}) = {\bf
Z}/2{\bf Z}\\

 [\sigma_{q}] &  \longmapsto &   \;\;\;\;\;\;\;\;\;\;
\mbox{lk}_{2}(p,q). \\
\end{array} $$
Since $\sigma_q|_{Y_p} = id_{Y_p} \Leftrightarrow \sigma_q(\sqrt{p})=\sqrt{p} \Leftrightarrow p \; \mbox{is a quadratic residue mod} \; q$, we have
$$ \displaystyle{(-1)^{\mbox{lk}_{2}(p,q)} = \left(
 \frac{p}{q} \right)}. \leqno{(2.1)}$$

As in the case of a link, the Legendre symbol is also interpreted as an ``intersection number". Firstly, we regard a prime $(q)$ as the dual of a ``meridian" around $(q)$, namely the Kummer character $\chi_{q} : \mbox{Gal}({\bf
Q}_{q}(\sqrt{q})/{\bf Q}_{q}) \rightarrow {\bf F}_{2}$ defined by
$\tau(\sqrt{q})/\sqrt{q} = (-1)^{\chi_{q}(\tau)}$. For $X_q := {\rm Spec}({\bf Z}[\frac{1}{q}])$, consider the
boundary map $ \partial \, : \,  H^{1}(X_{q},{\bf F}_{2})
\stackrel{\sim}{\rightarrow} H^{2}_{q}(\mbox{Spec}({\bf
Z}_{q}),{\bf F}_{2}) \simeq \mbox{Hom}(\mbox{Gal}({\bf Q}_{q}(\sqrt{q})/{\bf Q}_{q}),{\bf F}_{2})$ and we define the class of a ``Seifert surface" of $(q)$ by the class
$[\Sigma] \in H^{1}(X_{q},{\bf F}_{2})$ such that $\partial ([\Sigma]) =
\chi_{q}.$ On the other hand, we identify a ``knot" $(p)$ with $p \in {\bf
Q}_{q}^{\times}/({\bf Q}_{q}^{\times})^{2} = H^{1}({\bf Q}_{q},{\bf
F}_{2})$ and define the class $[(p)] \in H^{2}_{c}(X_{q},{\bf F}_{2})$ by its image under the map $H^{1}({\bf Q}_{q},{\bf F}_{2}) \rightarrow
H^{2}_{c}(X_{q},{\bf F}_{2})$, where $H_c^*$ stands for the compactly-supported \'{e}tale cohomology taking the infinite prime into account ([47,II]).
Following the case of a link, we define the mod 2 linking number by
the cup product $[(p)]\cup [\Sigma]$:
$$ \begin{array}{ccc}
H^{2}_{c}(X_{q},{\bf F}_{2}) \times H^{1}(X_{q},{\bf F}_{2}) &
\stackrel{\cup}{\longrightarrow}  & H^{3}_{c}(X_{q},{\bf F}_{2}) =  {\bf
Z}/2{\bf Z} \\
([(p)],[\Sigma])  & \longmapsto &  [(p)]\cup [\Sigma] = \mbox{lk}_{2}(p,q).
\\
\end{array} $$
By the construction, we obtain again the relation (2.1) (cf. [87]). Hence, the Legendre symbol is nothing but the mod 2 linking number and the Gauss reciprocity law corresponds to the symmetry of the linking number:
$$ \begin{array}{ccc}
\mbox{linking number} & \longleftrightarrow & \mbox{Legendre symbol},\\
{\rm lk}(K,L) = {\rm lk}(L,K) & \longleftrightarrow & \displaystyle{\left(\frac{p}{q}\right) = \left(\frac{q}{p}\right)} \;\; (p,q \equiv 1 \; \mbox{mod}\; 4).
\end{array}
\leqno{(2.2)}
$$

Among many proofs Gauss gave to the quadratic reciprocity law, there is a way to express the Legendre symbol using a Gaussian sum:
$$ \left( \sum_{x \in {\bf F}_p} \zeta^{x^2}\right)^{q-1} = \left( \frac{p}{q} \right) \leqno{(2.3)}$$
where $\zeta$ is a primitive $p$-th root of unity in $\overline{\bf F}_q$.
 Note that the Gaussian sum $\sum_{x \in {\bf F}_p}\zeta^{x^2}$ is an analogue over ${\bf F}_p$ of the Gaussian integral $\int_{\bf R}e^{-x^2}dx. $
On the other hand, Gauss showed the following integral expression of the linking number in his study of electro-magnetic theory ([14]):  
$$ 
\displaystyle{\int_{x \in K} \int_{y \in L} \omega(x-y)} = {\rm lk}(K,L) 
$$
where $\omega= (4\pi||x||^{3})^{-1}(x_1dx_2\wedge dx_3+x_2dx_3\wedge dx_1+x_3dx_1\wedge dx_2).$  We now rewrite Gauss' integral formula in a gauge-theoretic way ([34, 3.3]). Namely, for a framed link $K_{1} \cup K_{2}$, we have 
$$\begin{array}{l}
\displaystyle{\int_{A({\bf R}^{3})}\exp
\left(\frac{\sqrt{-1}}{4\pi}\int_{{\bf R}^{3}} a\wedge da +
\sqrt{-1}\int_{K_{1}} a +  \sqrt{-1}\int_{K_{2}} a\right){\cal D}a} \\
\;\;\;\;\;\;\;\;\;\;\;\;\;\;\;\;\;\;\;\; = \displaystyle{\exp\left( \pi\sqrt{-1} \sum_{1\leq i,j\leq 2}{\rm lk}(K_i,K_j)\right)}
\end{array}
\leqno{(2.4)}$$
Here the integral of the l.h.s is the Feynman integral over the space $A({\bf R}^{3})$
 of ${\bf R}$-valued differential 1-forms on ${\bf R}^{3}$. Since the integrals $\int_{{\bf R}^{3}}a\wedge da$ and $\int_{K_i}a$ define 
quadratic and linear forms on $A({\bf R}^{3})$ respectively, the
Feynman integral in (2.4) is regarded as an infinite dimensional analogue 
of the Gaussian integral. Thus we may observe an analogy between (2.3) and (2.4).\\
\\
{\bf 3. Link groups and pro-$l$ Galois groups with restricted ramification}
\\

As explained in Section 1, our basic idea is to regard the Galois group $G_{S} = \pi_1({\rm Spec}({\bf Z})\setminus S))$, $S=\{(p_1),\cdots ,(p_n)\}$, as an analogue of a link group $G_L=\pi_1(S^3 \setminus L)$, $L=K_1\cup \cdots \cup K_n$. Since the group $G_S$ is so huge, we consider its maximal pro-$l$ quotient $G_S(l)$ for a prime number $l$. Due to the works by I. \v{S}afarevi\v{c}, H. Koch etc, one has better understanding on pro-$l$ extensions of number fields ([33]). In fact, we find an analogy between a theorem by Koch on $G_S(l)$ and a theorem by J. Milnor on $G_L$. In the following, we denote by $\{G^{(d)}\}_{d\geq 1}$ the lower central series of a topological group $G$, which is defined by $G^{(1)}:=G, G^{(d+1)}:=[G^{(d)},G]$ = the closed subgroup topologically generated by $[a,b], a \in G^{(d)}, b \in G$.
\\

Let  $L = K_{1}\cup \cdots \cup K_{n}$ be  a link of $n$ components in
$S^{3}$ and let $G_{L}$ be the link group $\pi_{1}(S^{3} \setminus L)$.
  After the
work of K.T. Chen, Milnor derived the following information on the
presentation of the nilpotent quotients of $G_{L}$.\\
\\
{\bf Theorem 3.1} ([48]). {\em Let $F$ be the free group generated by $x_1,\dots ,x_n$. For each  $d \geq 1$, there is a word
$y_{i}^{(d)}$ of $x_{1},\cdots , x_{n}$ such that
$$ \begin{array}{l}
y_i^{(d)} \equiv y_i^{(d+1)} \; \mbox{mod}\; F^{(d)}\; (1\leq i\leq n),\\
G_{L}/G_{L}^{(d)} = \langle x_{1},\cdots ,x_{n} \; | \; [x_{i},
y_{i}^{(d)}] = 1  (1 \leq i \leq n), \; F^{(d)} = 1 \rangle. 
\end{array}
$$
Here $x_{1},\cdots , x_{n}$ is the basis of the free group $F$ where
$x_{i}$ represents a meridian $\alpha_{i}$ of $K_{i}$ and $y_{i}^{(d)}$
 is the word in $F$ representing the image of a longitude $\beta_{i}$ of
$K_{i}$ in $G_{L}/G_{L}^{(d)}$.}
{\em We also have the following relation}:
$$ \beta_{j} \equiv \prod_{i \neq j}\alpha_{i}^{\mbox{lk}(K_{i},K_{j})}
\; \mbox{mod} \; G_{L}^{(2)}.$$
\\
Milnor's Theorem 3.1 is extended for a link in any homology 3-sphere ([84]). The following theorem tells us that any link in a homology 3-sphere looks like a pure braid link after the pro-$l$ completion of the link group where $l$ is a prime number.\\
\\
{\bf Theorem 3.2} ([21]). {\em Let $L = K_1\cup\cdots\cup K_n$ be a link of $n$-components in a homology $3$-sphere $M$. For a prime number $l$,  let $\widehat{G_{L}}(M)$ denote the pro-$l$ completion of the link group $G_{L}(M) = \pi_{1}(M \setminus L)$. Then the pro-$l$ group $\widehat{G_{L}}(M)$ has the following presentation}:
$$ \widehat{G_{L}}(M) = \langle x_{1}, \cdots , x_{n} | [x_{1},y_{1}] = \cdots = [x_{n},y_{n}] = 1 \rangle.$$
where  $x_{i}$ represents a meridian of  $K_{i}$ and $y_i$ is the pro-$l$ word  representaing a longitude of $K_{i}$. \\

Applying the method by D. Anick ([1]) to Theorem 3.2,  we obtain the following theorem (pro-$l$ version of a conjecture of K. Murasugi ([41],[43])) on the structure of the lower central series of the link group $\widehat{G_{L}}$ ($M = S^{3}$). Let $D_L(l)$ be the mod $l$ linking diagram, namely the graph whose vertices are components of $L$, with two vertices $K_i$ and $K_j$ being joined by an edge if ${\rm lk}(K_{i},K_{j}) \not\equiv 0$ mod $l$.\\
\\
{\bf Theorem 3.3} ([21]). {\em If the graph $D_{L}(l)$ is connected, $ \widehat{G_{L}}^{(d)}/\widehat{G_{L}}^{(d+1)}$ is a free ${\bf Z}_{l}$-module whose ${\bf Z}_{l}$-rank $a_{d}$ is given by}:
$$ \prod_{d \geq 1}(1-t^{d})^{a_{d}} = (1-t)(1-(n-1)t).$$
{\em In particular, we have an isomorphism  $\widehat{G_{L}}^{(d)}/\widehat{G_{L}}^{(d+1)}  \simeq \widehat{F}_{1}^{(d)}/\widehat{F}_{1}^{(d+1)} \times \widehat{F}_{n-1}^{(d)}/\widehat{F}_{n-1}^{(d+1)}$ for $d\geq 1$ where $\widehat{F}_{r}$ stands for a free pro-$l$ group of rank $r$.}\\

 On the other hand, for a given prime number $l$, let $S = \{ (p_{1}),\cdots ,(p_{n}) \}$ be a finite set of $n$ distinct primes such that $p_i \equiv 1$ mod $l$ $(1\leq i\leq n)$. Let $G_{S}(l)$ denote the maximal pro-$l$ quotient of $\pi_{1}(\mbox{Spec}({\bf Z}) \setminus S)$, namely the Galois group $\mbox{Gal}({\bf Q}_{S}(l)/{\bf Q})$ of the maximal pro-$l$
 extension ${\bf Q}_{S}(l)$ over ${\bf Q}$ unramified outside $l$.  Koch derived the following information on the presentation of $G_{S}(l)$.\\
\\
{\bf Theorem 3.4} ([33], [51]). {\em The pro-$l$ group $G_S(l)$ has the following presentation}:
$$ G_{S}(l) = \langle x_{1},\cdots , x_{n} \; | \;
x_{i}^{p_{i}-1}[x_{i},y_{i}] = 1 (1\leq i \leq n) \rangle,$$
{\em where $x_{i}$ represents a monodromy $\tau_{i}$ over $p_{i}$ and
$y_{i}$ represents an extension of the Frobenius automorphism
$\sigma_{i}$ over $p_{i}$.}
{\em The mod $l$ linking number ${\rm lk}_l(p_i,p_j) \in {\bf F}_l$ is given by} 
$$ \sigma_{j} \equiv \prod_{i \neq j}\tau_{i}^{{\rm lk}_{l}(p_{i},p_{j})} \; \mbox{mod} \;
G_{S}(l)^{(2)},\;\; \zeta_{l} ^{{\rm lk}_{l}(p_{i},p_{j})} = \left( \frac{p_{j}}{p_{i}}
\right)_{l}$$
{\em where $\zeta_{l}$ a suitably chosen primitive $l$-th root of $1$
 and $(p_{j}/p_{i})_{l}$ is the $l$-th power residue symbol in
${\bf Q}_{p_{i}}$.}\\
\\
The analogy between Theorem 3.1 (or 3.2) and Theorem 3.4 is clear and explains the analogy between the linking number and the power residue symbol group-theoretically. We observe that the relations $[x_i,y_i]=1$ and $x_{i}^{p_{i}-1}[x_{i},y_{i}] = 1$ in the presentations of $G_L$ and $G_S(l)$ are  coming from those (1.2) of the local fundamental groups  $\pi_1(\partial V_{K_i})$ ($V_{K_i}$ being a tubular neighborhood of $K_i$) and $\pi_{1}^{tame}({\rm Spec}({\bf Q}_{p_{i}}))$ respectively. Furthermore, as Milnor's theorem is extended for a link in any homology 3-sphere, Koch's theorem can be generalized for a finite set of primes in a number field: Let $k$ be a number field containing a primitive $l$-th root of unity and $S = \{ \frak{p}_{1},\cdots ,\frak{p}_{n}\}$  a finite set of $n$ distinct primes such that $N\frak{p}_{i} := \#({\cal O}_{k}/\frak{p}_{i}) \equiv 1$ mod $l$ ($1\leq i \leq n$). Set $B_{S} := \{ \alpha \in k^{\times} | (\alpha ) = \frak{a}^{l} \,(\exists \frak{a} : \,\mbox{ideal of}\; k),\, \alpha \in (k_{\frak{p}_{i}}^{\times})^{l} \; (1\leq i\leq n) \}/(k^{\times})^{l}$.
\\
\\
{\bf Theorem 3.5} ([33]). {\em Assume that $k$ contains a primitive $l$-th root of unity and $B_{S} = 1$ and that the class number of $k$ is prime to $l$. Then the maximal pro-$l$ quotient $G_{S}(k)(l)$  of  $\pi_{1}({\rm Spec}({\cal O}_{k}) \setminus S)$ has the following presentation}:
$$G_{S}(k)(l) = \langle x_{1},\cdots , x_{n} \; | \; x_{1}^{N\frak{p}_{1}-1}[x_{1},y_{1}] = \cdots = x_{n}^{N\frak{p}_{n}-1}[x_{n},y_{n}] = 1  \rangle,$$
 {\em where  $x_{i}$ represents a monodromy $\tau_i$ over  $\frak{p}_{i}$ and  $y_{i}$ represents an extension of the Frobenius automorphism
$\sigma_{i}$ over $\frak{p}_{i}$.} \\

 An arithmetic analogue of the Murasugi conjecture for the Galois group $G_{S}(l)$ is a more delicate problem. In fact, we can not expect an analogy for primes in the case $l > 2$, since the linking number ${\rm lk}_{l}(p_{i},p_{j})$ is not symmetric. However, J. Labute showed recently the following theorem on the $l$-lower central series $\{ G_S(l)_d\}$ of $G_{S}(l)$ which is defined by $G_{S}(l)_{1} := G_{S}, G_{S}(l)_{d+1} := (G_{S}(l)_{d})^{l}[G_{S}(l)_{d},G_{S}]$. The mod $l$ linking graph $D_{S}(l)$ of $S$ is defined by the graph whose vertices are primes of $S$, with two vertices $(p_i)$ and $(p_j)$ being joined by an edge if ${\rm lk}_{l}(p_{i},p_{j}) \neq 0$.\\
\\
{\bf Theorem 3.6} ([37]). {\em Suppose that $l > 2$ and set $a_{d} = \dim_{{\bf F}_{l}}G_{S}(l)_{d}/G_{S}(l)_{d+1}$ for $d \geq 1$. Assume the following conditions on the graph} $D_{S}(l)$: {\em $(1)$ The vertices $p_{1},p_{2},\cdots ,p_{n}$ form a circuit  $p_{1}p_{2}\cdots p_{n}p_{1}$. $(2)$ If $i,j$ are both odd, $p_{i}p_{j}$ is not a edge of $D_{S}(p)$. $(3)$ $l_{12}l_{23}\cdots l_{n-1,n}l_{n,1} \neq l_{1n}l_{2,1}\cdots l_{n,n-1}$ for $l_{ij} := {\rm lk}_{l}(p_{i},p_{j})$. Then we have the following formula}:
$$ \prod_{d\geq 1}(1-t^{d})^{a_{d}} = (1-t)(1-nt+nt^{2}).$$
\\
For the case that $p=2$ and $p_{i} \equiv 1$ mod 4, the mod 2 linking number ${\rm lk}_{2}(p_{i},p_{j})$ is symmetric and we may expect an analogue of the Murasugi conjecture for the Zassenhaus filtration of $G_S(2)$. However we do not have such a result in terms of the graph $D_{S}(2)$ yet. For the study of $G_S(l)$ by the linking graph $D_S(l)$, we also refer to [74].\\
\\
{\bf 4. Milnor invariants and multiple power residue symbols}
\\

The Milnor $\overline{\mu}$-invariants for a link were introduced by Milnor ([48]) and the relation with  Massey products was
 proven by V. Turaev ([84]) and R. Porter ([67]). Based on the analogy between a link group and a pro-$l$ Galois group with restricted ramification in Section 3, we can introduce arithmetic analogues of the Milnor invariants and Massey
products for prime numbers  which generalize the power residue symbols and R\'{e}dei triple symbols ([69]). The basic tool is a pro-$l$ version ([23]) of the Fox free differential calculus ([5]).\\

 Let $L = K_{1}\cup \cdots \cup K_{n}$ be a link of $n$ component in $S^3$ (or a homology 3-sphere), and let $F$ be the free group generated by the words $x_1,\dots ,x_n$ where each $x_i$ represents a meridian
$\alpha_{i}$ of $K_i$ as in Section 3. Let ${\bf Z}\langle \langle X_{1},\cdots ,X_{n} \rangle \rangle$ be the
ring of formal power series in non-commuting variables $X_{1}, \cdots ,
X_{n}$ and let $M : F \hookrightarrow {\bf Z}\langle \langle X_{1},\cdots ,
X_{n} \rangle \rangle^{\times}$ be the Magnus embedding defined by sending $x_{i}$ to $1 + X_{i}$. We then expand the
longitude $y_j^{(d)}$ of $K_j$ by meridians:
$$ M(y_{j}^{(d)}) = 1 + \sum \mu(i_{1}\cdots i_{r} j)X_{i_{1}}\cdots
X_{i_{r}}.$$
Alternative description of $\mu(i_{1}\cdots i_{r}j)$ is given using the
Fox free differential calculus on the group ring ${\bf Z}[F]$ ([5]):
$$ \mu(i_{1}\cdots i_{r}j) = \epsilon \left( \frac{\partial^{r}
y_{j}^{(d)}}{\partial x_{i_{1}}\cdots \partial x_{i_{r}}} \right) $$
where $\epsilon$ is the augmentation homomorphism ${\bf Z}[F]
\rightarrow {\bf Z}$. For a multi-index $I$, let $\Delta(I)$ be the ideal of ${\bf Z}$ generated by $\mu(J)$ where $J$
 ranges over all cyclic permutations of proper subsequences of $I$, and
we define the  Milnor $\overline{\mu}$-invariant by
$$ \overline{\mu}(I) := \mu(I) \; \mbox{mod} \; \Delta(I).$$
 Milnor showed the followings.
\\
\\
{\bf Theorem 4.1} ([48]). (1) $\overline{\mu}(ij) =
\mbox{lk}(K_{i},K_{j})$.\\
(2) {\em The Milnor $\overline{\mu}$-invarints $\overline{\mu}(I)$ are
isotropy invariants of $L$ for $|I| < d$, where $|I|$ stands for the
length of $I$.}\\
(3) {\em We have the following symmetry relations}:\\
({\em cyclic symmetry})$\;\;$  $\overline{\mu}(i_{1}\cdots i_{r}) =
\overline{\mu}(i_{2}\cdots i_{r}i_{1}) = \cdots =
\overline{\mu}(i_{r}i_{1}\cdots i_{r-1}).$\\
({\em shuffle relation}) {\em If $Sh$ denotes the set of all proper
shuffles of multi-indices $I$ and $J$ with $|I|, |J| \geq 1$, then }
$$ \sum_{H \in Sh}\overline{\mu}(Hk) \equiv 0 \;\;  \mbox{mod} \;\;
\mbox{g.c.d}\{\Delta(Hk) | H \in Sh \}.$$
\\
{\bf Example 4.2.} (1) Let $L=K_1\cup K_2$ be the  Whitehead link. Then we have $\overline{\mu}(I)=0$ for $|I|=2,3$ or $|I|\geq 5$, and $\overline{\mu}(1122)=1, \overline{\mu}(1212)=-2$.\\
(2) Let $L= K_{1}\cup K_{2}\cup K_{3}$ be the  Borromean rings. Then
 we have $\mu(ij) = {\rm lk}(K_{i},K_{j}) = 0$ for
all $i, j$ and $\mu(ijk) = \pm 1$ for any permutation $ijk$ of $123$. \\
\vspace{2.8cm}
$$ \mbox{Whitehead link} \;\;\;\;\;\;\;\;\;\;\;\;\;\;\;\;\;\;\;\;\;\;\;\;\;\;\;\;  \mbox{Borromean ring}$$
\\
 
As the linking number is an invariant of an abelian covering, the Milnor invariants are interpreted as covering invariants in
 nilpotent coverings over $S^{3}$ as follows. Let
$N_{r}(R)$ be the group of all upper triangular $r \times r$ matrices
with 1 along the diagonal entries over a ring $R$. For a multi-index $I
= (i_{1}\cdots i_{r})$ ($r \geq 2$), we define a nilpotent
representation $\rho_{I} \, : \, F \longrightarrow N_{r}({\bf
Z}/\Delta(I))$ by
$$ \rho_{I}(f) := \left( \begin{array}{ccccc}
1 & \epsilon\left(\frac{\partial f}{\partial x_{i_{1}}}\right) &
\epsilon\left(\frac{\partial^{2}f}{\partial x_{i_{1}}\partial x_{i_{2}}}\right)  &
\cdots & \epsilon\left(\frac{\partial^{r-1}f}{\partial x_{i_{1}}\cdots
\partial x_{i_{r-1}}}\right) \\
  &  1 & \epsilon\left(\frac{\partial f}{\partial x_{i_{2}}}\right) & \cdots &
\epsilon\left(\frac{\partial^{r-2}f}{\partial x_{i_{2}} \cdots \partial
x_{i_{r-1}}}\right) \\
 & & \ddots & \ddots & \vdots \\
 & \mbox{\Huge 0}&        & 1 & \epsilon\left(\frac{\partial f}{\partial
x_{i_{r-1}}}\right) \\
 & &        &  & 1  \end{array} \right)  \; \; {\rm mod} \;\;
\Delta(I).$$
\\
{\bf Theorem 4.3.} (cf. [63]) {\em Notations being as above,}\\
(1) {\em the representation $\rho_{I}$ factors through the link group $G_{L}$, and
further it gives a surjective representation of $G_{L}$ onto $N_{r}({\bf
Z}/\Delta(I))$ if $i_{1}, \cdots , i_{r-1}$ are distinct each other.}\\
(2) {\em Suppose $i_{1}, \cdots , i_{r-1}$ are distinct each other. If
$X_{r}$ denotes the covering of $X_{L} = S^{3} \setminus L$
 corresponding to ${\rm Ker}(\rho_{I})$, then  $X_{r} \rightarrow X_{L}$
 is a Galois covering ramified over $K_{i_{1}}\cup \cdots \cup
K_{i_{r-1}}$ with Galois group $N_{r}({\bf Z}/\Delta(I))$, and for a
longitude $\beta_{i_{r}}$ around $K_{i_{r}}$, we have}
$$\rho_{I}(\beta_{i_{r}}) = \left( \begin{array}{ccccc}
1 & 0 & 0  & \cdots & \overline{\mu}(I) \\
  &  1 & 0 & \cdots & 0  \\
 & & \ddots & \ddots & \vdots \\
 & \mbox{\Huge 0}&        & 1 & 0  \\
 & &        &  & 1  \end{array} \right).$$
\\

The Milnor invariants are also interpreted in terms of the Massey products in the cohomology of the complement $X_L = S^3 \setminus L$. Let $\xi_{1},\dots ,\xi_{n}$
 be the basis of $H^{1}(X_{L},{\bf Z})$ dual to the meridians
$\alpha_{1},\dots ,\alpha_n$, and let $\eta_{j} \in H_{2}(X_{L},{\bf Z})$ be the class
realized by the boundary $\partial V_{K_{j}}$, $V_{K_{j}}$ being a tubular neighborhood around $K_{j}$. Then we have the following\\
\\
{\bf Theorem 4.4} ([67],[84]). {\em For $I = (i_{1}\cdots
i_{r})$,  there is a defining  system
$M$ for the Massey product $\langle \xi_{i_{1}}, \cdots , \xi_{i_{r}}
\rangle$ in $H^{2}(X_{L},{\bf Z}/\Delta(I))$ so that we have}
$$ \langle \xi_{i_{1}}, \cdots , \xi_{i_{r}} \rangle_{M} (\eta_{j}) =
\left\{ \begin{array}{rl}
(-1)^{r} \overline{\mu}(I), & \quad \mbox{for $j = i_{r} \neq i_{1}$}\\
(-1)^{r+1} \overline{\mu}(I), & \quad \mbox{for $j = i_{1} \neq
i_{r}$}\\
0, & \quad \mbox{otherwise.} \end{array}\right.  $$
\vspace{.05cm}

On the other hand, let $S= \{(p_1),\dots ,(p_n)\}$, $p_i \equiv 1$ mod $l$, be a set of $n$ distinct primes as in Section 3, and let $\hat{F}$ be the
free pro-$l$ group generated by the words $x_{1},\dots ,x_n$ where each $x_i$ represents the monodromy $\tau_{i}$ over $p_i$ in the Galois group $G_{S}(l)$.
 We set $e_S:= \max\{ e \,|\, p_i \equiv 1 \; \mbox{mod}\; l^e \, (\forall i)\}$ and fix $m= l^e$ $(1\leq e\leq e_S)$. 
Let ${\bf Z}_l\langle \langle
X_{1},\cdots , X_{n} \rangle \rangle$ be the ring of formal power series in non-commuting
variables $X_{1}, \cdots , X_{n}$ over ${\bf Z}_l$ and let $\hat{M} : \hat{F} \hookrightarrow {\bf Z}_l\langle \langle
X_{1},\cdots , X_{n} \rangle \rangle^{\times}$ be the Magnus embedding
 over ${\bf Z}_l$ defined by sending $x_i$ to $1+X_i$. We then
expand the Frobenius $y_j$ over $p_j$ by monodromies:
$$ \hat{M}(y_{j}) = 1 + \sum \hat{\mu}(i_{1}\cdots i_{r}j)X_{i_{1}}\cdots
X_{i_{r}}$$
and define
$$ \mu_m(I) := \hat{\mu}(I) \; \mbox{mod}\; m.$$
Alternative description of $\hat{\mu}(i_{1}\cdots i_{r}j)$ is given  
using the pro-$l$ Fox free differential calculus ([23]) on the
completed group ring ${\bf Z}_{l}[[\hat{F}]]$:
$$ \mu_{m}(i_{1}\cdots i_{r}j) = \epsilon_{l} \left( \frac{\partial^{r}
y_{j}}{\partial x_{i_{1}}\cdots \partial x_{i_{r}}} \right) $$
where $\epsilon_{l}$ is the augmentation homomorphism ${\bf
Z}_{l}[[\hat{F}]] \rightarrow {\bf Z}_{l}.$ For a multi-index $I$, $1 \leq |I| \leq l^{e_{S}}$, let $\Delta(I)$
 be the ideal of ${\bf Z}/m{\bf Z}$ generated by the binomial
coefficients $l^{e_{S}}\choose t$ and $\mu_{m}(J)$ where $1 \leq t \leq
|I|$ and $J$ ranges over all cyclic permutations of proper subsequences
of $I$. We then define the arithmetic Milnor $\overline{\mu_{m}}$-invariant
by
$$ \overline{\mu_{m}}(I) = \mu_{m}(I) \; \; \mbox{mod} \; \Delta(I).$$
Our theorem is stated in the following form analogous to Theorem 4.1.
\\
\\
{\bf Theorem 4.5}  ([50],[51],[54]). (1) {\em We have
$\zeta_{m}^{\mu_{m}(ij)} = \displaystyle{\left( \frac{p_{j}}{p_{i}}
\right)_{m}}$ where $\zeta_{m}$ is a primitive root of $1$ in ${\bf
Q}_{p_{i}}$}.\\
 (2) {\em Let $I$ be a multi-index with $2 \leq |I| \leq l^{e_{S}}$.
Then  $\overline{\mu_{m}}(I)$ are invariants determined by the Galois
group $G_{S}(l)$.}\\
(3) {\em  Let $r$ be an integer with $2 \leq r \leq l^{e_{S}}$. }
 {\it For multi-indices $I = (i_{1}\cdots i_{s})$ and $J =
(j_{1}\cdots j_{t})$ with $s + t = r - 1, s,t \geq 1$, 
 $Sh$ denotes the set of all proper shuffles of $I$ and $J$.  Then we have}:
$$ \sum_{H \in Sh}\overline{\mu_{m}}(Hk) \equiv 0 \;\;  \mbox{mod} \;\;
\mbox{g.c.d}\{\Delta(Hk) | H \in Sh \}.$$
\vspace{.05cm}

As in the case of a link, the Milnor $\overline{\mu_m}$-invariants describes the decomposition law of a prime in
  nilpotent extensions of ${\bf Q}$.  For a multi-index $I = (i_{1}\cdots
i_{r})$, $2 \leq r < l^{e_{S}}$ such that $\Delta(I) \neq {\bf Z}/m{\bf
Z}$, we define a representation $\rho_{m,I} \, : \, \hat{F} \longrightarrow N_{r}(({\bf Z}/m{\bf
Z})/\Delta(I))$ by
$$ \rho_{m,I}(f) := \left( \begin{array}{ccccc}
1 & \epsilon\left(\frac{\partial f}{\partial x_{i_{1}}}\right)_{m} &
\epsilon\left(\frac{\partial^{2}f}{\partial x_{i_{1}}\partial x_{i_{2}}}\right)_{m}
  & \cdots & \epsilon\left(\frac{\partial^{r-1}f}{\partial x_{i_{1}}\cdots
\partial x_{i_{r-1}}}\right)_{m} \\
  &  1 & \epsilon\left(\frac{\partial f}{\partial x_{i_{2}}}\right)_{m} & \cdots &
\epsilon\left(\frac{\partial^{r-1}f}{\partial x_{i_{2}} \cdots \partial
x_{i_{r-1}}}\right)_{m} \\
 & & \ddots & \ddots & \vdots \\
 & \mbox{\Huge 0}&        & 1 & \epsilon\left(\frac{\partial f}{\partial
x_{i_{r-1}}}\right)_{m} \\
 & &        &  & 1  \end{array} \right)  \; \; {\rm mod} \;\;
\Delta(I)$$
where we set $\epsilon(\alpha)_{m} = \epsilon_{l}(\alpha)$ mod $m$.\\
\\
{\bf Theorem 4.6} ([54]). {\em Notations being as above,}\\
(1) {\em the representation $\rho_{m,I}$ factors through $G_{S}(l)$, and
it gives a surjective representation of $G_{S}(l)$ onto $N_{r}(({\bf
Z}/m{\bf Z})/\Delta(I))$ if $i_{1}, \cdots , i_{r-1}$ are distinct each
other.}\\
(2) {\em Suppose $i_{1}, \cdots , i_{r-1}$ are distinct each other. If
$k_{r}$ denotes the extension of ${\bf Q}$ corresponding to ${\rm
Ker}(\rho_{m,I})$, then  $k_{r}/{\bf Q}$ is a Galois extension ramified
over $p_{i_{1}}, \cdots , p_{i_{r-1}}$ with Galois group $N_{r}(({\bf
Z}/m{\bf Z})/\Delta(I))$, and for the Frobenius automorphism
$\sigma_{i_{r}}$ over $p_{i_{r}}$, we have}
$$\rho_{I}(\sigma_{i_{r}}) = \left( \begin{array}{ccccc}
1 & 0 & 0  & \cdots & \overline{\mu_{m}}(I) \\
  &  1 & 0 & \cdots & 0  \\
 & & \ddots & \ddots & \vdots \\
 & \mbox{\Huge 0}&        & 1 & 0  \\
 & &        &  & 1  \end{array} \right).$$
\\
{\bf Example 4.7} ([50],[51],[54],[86]). Suppose that $l =2$ and  $p_{1}, p_{2}, p_{3}$ are distinct prime numbers congruent to 1  mod 4. By Theorem 4.6, there  is a Galois extension $k$ of ${\bf Q}$ ramified over $p_1$ and $p_2$ with Galois group $N_3({\bf F}_2)$ = the dihedral group of order 8. Here noting that $\mu_2(123)=0 \Leftrightarrow \sigma_3|_{k} = id \Leftrightarrow ``p_3 \; \mbox{is completly decomposed in}\;  k/{\bf Q}"$, we see that the  triple symbol $[p_1,p_2,p_3]$ introduced by L. R\'{e}dei ([69]) is interpreted as our Milnor invariant:
$$ [p_{1},p_{2},p_{3} ] = (-1)^{\mu_{2}(123)}.$$
For example, for $S = \{13,61,937\}$, we have the same relations which the Borromean ring satisfies (Example 4.2 (2)): $\mu_{2}(ij) = {\rm lk}_{2}(p_{i},p_{j}) = 0$ for any $i\neq j$,
$\mu_{2}(ijk) = 1$ for any permutation of $123$ and $\mu_2(ijk)=0$ otherwise. This triple may be called mod 2  ``Borromean primes". \\

As in the case of a link, our arithmetic Milnor invariants
are interpreted as Massey products in the \'{e}tale cohomology of the complement $X_S := {\rm Spec}({\bf Z}) \setminus S$. Let $\xi_{1},\dots, \xi_{n}$ be the basis of
$H^{1}(X_{S},{\bf Z}/m{\bf Z})$ dual to the monodromies $\tau_1,\dots ,\tau_n$ and let $\eta_j \in H_2(X_S,{\bf Z}/m{\bf Z})$ be the image of the canonical generator of $H_2({\rm Spec}({\bf Q}_{p_i}),{\bf Z}/m{\bf Z})$.\\
\\
{\bf Theorem 4.9} ([54]). {\em  Let $I = (i_{1}\cdots i_{r})$ be a
multi-index with $2 \leq r \leq l^{e_{S}}$. Then there is a defining system $M$ for the Massey product $\langle \xi_{i_{1}},\cdots ,
\xi_{i_{r}} \rangle$ in $H^{2}(X_S,({\bf Z}/m{\bf Z})/\Delta (I))$
 so that we have}
$$\langle \xi_{i_{1}}, \cdots , \xi_{i_{r}} \rangle_{M}(\eta_{j}) =
\left\{ \begin{array}{rl}
(-1)^{r} \overline{\mu_{m}}(I), & \quad \mbox{for $j = i_{r} \neq
i_{1}$}\\
(-1)^{r+1} \overline{\mu_{m}}(I), & \quad \mbox{for $j = i_{1} \neq
i_{r}$}  \\
0, & \quad \mbox{otherwise.} \end{array}\right.$$
\\
This generalizes the well-known relationship between power residue symbols
and  cup products ([33, 8.11]) and in particular gives a cohomological
interpretation of the R\'{e}dei triple symbol. Though the Milnor invariants $\overline{\mu_m}(I)$ are defined by using the Galois group $G_S(l)$, the Massey product interpretation shows that they depend only on $X_S$ (and $l$). Finally, we  give an extension of the Milnor invariants to number fields. Let $k$ and $S$ be a pair of a number field and a finite set of primes which satisfies the conditions in Theorem 3.5. By Theorem 3.5, we can introduce the Milnor $\overline{\mu_m}$-invariants in the same manner as in the case of ${\bf Q}$. We then define the multiple power residue symbols for primes by
$$ [\frak{p}_{i_1},\dots ,\frak{p}_{i_r}] := \zeta_l^{\mu_m(i_1\cdots i_r)}$$
which may be regarded as a multiple generalization of the Legendre and R\'{e}dei symbols to a number field.\\
\\
{\bf 5. Homology groups and ideal class groups}\\

The notion of an ideal class group goes back to Gauss' theory on binary quadratic forms. In terms of number fields, Gauss showed that ideal classes of a quadratic number field form an finite abelian group and determined the structure of genera classifying ideal classes. These notions and results are clearly understood and generalized naturally from the standpoint of the analogy with knot theory. In the rest of this paper, chain and homology groups are assumed to have the integral coefficient ${\bf Z}$, and denoted simply by $C_{*}(M), H_{*}(M)$ etc.\\

For a 3-manifold $M$, all knots in $M$ generate the group $Z_1(M)$ of 1-cycles and the 1st homology group $H_1(M)$ is defined to be the quotient of $Z_1(M)$ by the subgroup made of the boundaries $\partial \Sigma$ for 2-cycles $\Sigma$. Similarly, for a number field $k$, all prime ideals of ${\cal O}_k$ generate the ideal group $I_k$ and the (narrow) ideal class group $H_k$ is defined to be the quotient of $I_k$ by the subgroup made of principal ideals $(a)$ for (totally positive) $a \in k^{\times}$:
$$ \begin{array}{ccc}
 C_{2}(M) \rightarrow C_{1}(M) & \longleftrightarrow &k^{\times} \rightarrow I_{k} \\
\Sigma \mapsto \partial \Sigma & &  a \; \mapsto (a) \\
\mbox{1st homology group}\; H_{1}(M) & \longleftrightarrow & \mbox{ideal class group} \;H_{k}.
\end{array}
$$

First, let us recall Gauss' theory of genera. Let $k/{\bf Q}$ be a quadratic number field ramified over $n$ distinct odd prime numbers $p_1,\dots ,p_n$. Two ideals $I$ and $J$ are said to be in the same genus if $\left( \frac{NI}{p_{i}}\right) = \left( \frac{NJ}{p_{i}}\right)$ holds for $1 \leq i \leq n$. Since $I$ and $J$ are in the same genus if they are in the same class, we can classify $H_k$ by the genus relation. Then Gauss's genus theory asserts that ideal classes in the genus of the identity $[{\cal O}_k]$ form $H_k^{2}$, and that the correspondence $[I] \mapsto ( \left( \frac{NI}{p_{1}}\right), \cdots ,\left(\frac{NI}{p_{n}}\right))$ induces an isomorphism
$$ H_k/H_k^2 \simeq \{(u_{i}) \in \{\pm 1\}^{n}\, | \prod_{i=1}^{n}u_{i} = 1\}
\simeq ({\bf Z}/2{\bf Z})^{n-1}.$$ 
 According to the analogies (2.2) and (5.1), we have an analogue of Gauss' genus theory for a link  as follows.  Let $f : M \rightarrow S^{3}$ be a double covering of closed 3-manifolds branched along a link $K_{1}\cup \cdots \cup K_{n}$ of $n$ components. We say that 1-cycles $a, b$ are in the same genus if  ${\rm lk}(f_{*}(a),K_{i}) \equiv {\rm lk}(f_{*}(b),K_{i}) \,{\rm mod} \, 2 \, (1\leq i\leq n)$ holds for $1\leq i\leq n$. Since $a$ and $b$ are in the same genus if they are in the same homology class, we can classify $H_1(M)$ by the genus relation. Then we can show ([52],[71],[73]) that homology classes in the genus of the identity form  $2H_{1}(M)$, and the correspondence $[c] \mapsto ({\rm lk}(f_{*}(c),K_{1}) \, {\rm mod} \, 2, \cdots , {\rm lk}(f_{*}(c),K_{n})\, {\rm mod} \, 2)$ induces an isomorphism 
$$H_{1}(M)/2H_{1}(M) \simeq \{ (t_i) \in ({\bf Z}/2{\bf Z})^{n}\, | \, \sum_{i=1}^{n} t_i = 0 \} \simeq ({\bf Z}/2{\bf Z})^{n-1}.$$ 
Thanks to Gauss' genus theory, the 2-Sylow subgroup $H_k(2)$ of the narrow ideal class group of the quadratic field $k$ has the form 
$$H_{k}(2) = \bigoplus_{i=1}^{n-1}{\bf Z}/2^{a_{i}}{\bf Z}, \; \;a_{i}\geq 1.$$ 
Hence the 2-rank of $H_{k}$ is determined by the number of ramified prime numbers in $k/{\bf Q}$. Since Gauss' time, it has been a problem to determine the whole structure of $H_k(2)$, namely to express the $2^{d}$-rank of $H_k(2)$ for $d >1$ in terms of the quantities related to $p_{1},\cdots ,p_{n}$ ([88]). Among many works on this problem, R\'{e}dei ([69]) expressed the 4-rank by using a matrix whose entries are given by the Legendre symbols involving $p_{i}$'s. According to the analogy (2.2), we find that R\'{e}dei's matrix is nothing but the arithmetic analogue of the mod 2 linking matrix. Therefore it would be a natural generalization to express the $2^{d}$-rank by using a ``higher linking matrix" involving the Milnor numbers in Section 4. In fact, we can show such a formula for a link, and then, imitating the method for a link, we can obtain an higher order generalization of Gauss' genus theory.\\

Let us generalize the setting a little bit. Let $L= K_{1}\cup \cdots \cup K_{n}$ be a link of $n$ components in $S^{3}$. For a given prime number $l$, let $\psi : G_{L} \rightarrow \langle t \, | \, t^{l} = 1 \rangle$ be the homomorphism sending each meridian of $K_{i}$ to $t$, and let $M$ be the completion of the $l$-fold cyclic covering of $X_{L}$ corresponding to the kernel of $\psi$. We assume that  $M$ is a rational homology 3-sphere. As in the case that $l=2$ described above, the homomorphism  $H_{1}(M) \ni [c] \mapsto ({\rm lk}(f_{*}(c),K_{1}) \; {\rm mod} \, l, \cdots , {\rm lk}(f_{*}(c),K_{n})\; {\rm mod} \, l) \in {\bf F}_{l}^{n}$ induces the isomorphism  $H_{1}(M)/(t-1)H_{1}(M) \simeq {\bf F}_{l}^{n-1}.$ Therefore, the $l$-Sylow subgroup $H_{1}(M)(p)$ of $H_{1}(M)$ regarded as a module over the complete discrete valuation ring $\widehat{\cal O} := {\bf Z}_{l}[t]/(t^{l-1}+\cdots +t+1) = {\bf Z}_{l}[\zeta]$ ($\zeta := t \, {\rm mod}\, (t^{l-1}+\cdots +t+1))$ has $\frak{p}:= (\zeta -1)$-rank $n-1$. Hence we have
$$\displaystyle{H_{1}(M)(l) = \bigoplus_{i=1}^{n-1}\widehat{\cal O}/\frak{p}^{a_{i}}}, \;\; a_{i}\geq 1$$
as $\widehat{\cal O}$-module. Set $\pi := \zeta - 1$. Using the Milnor number $\mu(I)$ ($|I|<d+1$) in Section 4, we define the $d$-th linking matrix $T_{L}^{(d)} =(T_{L}^{(d)}(i,j))$  for $d\geq 2$ by 
$$ T_{L}^{(d)}(i,j)  := \left\{ 
\begin{array}{ll} \displaystyle{-\sum_{r=1}^{d-1} \sum_{{\scriptstyle 1\leq i_{1},\cdots ,i_{r}\leq n}\atop{\scriptstyle i_{r}\neq i}} \mu(i_{1}\cdots i_{r}i)\pi^{r}} & \quad i = j \\
 \mu(ji)\pi + \displaystyle{\sum_{r=1}^{d-2}\sum_{1\leq i_{1},\cdots ,i_{r}\leq n} \mu(i_{1}\cdots i_{r}ji)\pi^{r+1}} & \quad i \neq j.
\end{array}\right.   $$ 
\\
{\bf Theorem 5.2} ([21]). {\em For $d\geq 2$, $T_{L}^{(d)}$ gives a presentation matrix of the  $\widehat{\cal O}/\frak{p}^{d}$-module $H_1(M)(l)\oplus \widehat{\cal O}/\frak{p}^{d}$. Let $\varepsilon_{1}^{(d)},\cdots ,\varepsilon_{n-1}^{(d)}, \varepsilon_{n}^{(d)}=0$ $(\varepsilon_{i}^{(d)}|\varepsilon_{i+1}^{(d)})$ be the elementary divisors of $T_{L}^{(d)}$.  Then the $\frak{p}^{d}$-rank $e_{d} = \# \{ i \, | \, a_{i} \geq d \}$  $H_{1}(M)(l)$ is given by}:
$$ e_d = \#\{ i \,|\, \varepsilon_{i}^{(d)} \equiv 0 \, {\rm mod}\, \frak{p}^{d}\} - 1.$$
 \\
{\bf Corollary 5.3} ([21]).  (1) {\em We have} 
$$e_{2} = n - 1 -{\rm rank}_{{\bf F}_{l}}(T_{L}\;{\rm mod} \; l)$$
 {\em where $T_{L} = (T_{L}(i,j))$ is the linking matrix defined by} $T_{L}(i,i) = -\sum_{j\neq i}{\rm lk}(K_{i},K_{j})$, $T_{L}(i,j)= {\rm lk}(K_{i},K_{j})$ ($i\neq j$).\\
(2) {\em Suppose that $n =2$. Under the assumption that $e_{d} = 1$ $(d\geq 1)$, we have}
$$
e_{d+1} = 1 \Longleftrightarrow \left\{
\begin{array}{rl}
\displaystyle{\sum_{r=1}^{d}\sum_{i_{1},\cdots ,i_{r-1} = 1,2}\mu(i_{1}\cdots i_{r-1}21)\pi^{r} \; \equiv \; 0 \; {\rm mod}\; \pi^{d+1}},\\
\displaystyle{\sum_{r=1}^{d}\sum_{i_{1},\cdots ,i_{r-1} = 1,2}\mu(i_{1}\cdots i_{r-1}12)\pi^{r} \; \equiv \; 0 \; {\rm mod}\; \pi^{d+1}}. \\
\end{array}\right.
 $$
\\
The prototype of the matrix $T_{L}^{(d)}$ is a presentation matrix of the reduced Alexander module $A_L$ of $L$ introduced by L. Traldi ([83]). In fact,  the essential point of our proof of Theorem 5.2 is to complete $A_L$ to get a compact module $\hat{A}_L$ over the Iwasawa algebra ${\bf Z}_{l}[[T]]$ and to show that the ``universal higher linking matrix" defined by replacing $\pi$ by $T$ and  $d$ by $\infty$ in the definition of $T_{L}^{(d)}$ gives a presentation matrix of the completed Alexander module $\hat{A}_L$. \\
\\
{\bf Example 5.4} ([21]). Let $L = K_{1} \cup K_{2}$ be the Whitehead link (Example 4.2, (1)). Then we have $e_{3} = 1, e_{4} = 0$ and hence $H_{1}(M)(l) = \widehat{\cal O}/\frak{p}^{3}$. \\

On the other hand, let $S = \{ (p_{1}),\cdots ,(p_{n})\}$, $p_{i} \equiv 1$ mod $l$, be a finite set of primes as in Section 4. Let $\psi : G_{S}(l) \rightarrow \langle t \, | \, t^{l} = 1\rangle$ be the homomorphism sending each monodromy over $p_i$ to $t$, and let $k$ be the cyclic extension of ${\bf Q}$ of degree $l$ corresponding to the kernel of $\psi$. Let $\mu_{l}$ be the group of $l$-th roots of unity and fix an embedding $\mu_{l} \subset {\bf Q}_{p_{i}}$. Then the homomorphism $H_{k} \ni [I] \mapsto ( \left( \frac{NI}{p_{1}}\right)_{l}, \dots ,\left(\frac{NI}{p_{n}}\right)_{l} ) \in \mu_{l}^{n}$ induces the isomorphism $H_{k}/(t-1)H_{k} \simeq {\bf F}_{l}^{n-1}$ ([25]). Hence the $l$-Sylow subgroup $H_{k}(l)$ of the narrow ideal class group $H_k$ has $\frak{p}$-rank $n-1$  so that we have  
$$\displaystyle{H_{k}(l) = \bigoplus_{i=1}^{n-1}\widehat{\cal O}/\frak{p}^{a_{i}}}, \;\; a_{i}\geq 1$$
as $\widehat{\cal O}$-module. The $\frak{p}^{d}$-rank $e_{d}$ is also described by the arithmetic higher linking matrix $T_{S}^{(d)}$, where  $T_{S}^{(d)}$ is defined by replacing  $\mu(I)$ by $\widehat{\mu}(I)$ in the definition of $T_{L}^{(d)}$. Let $\varepsilon_1^{(d)},\cdots ,\varepsilon_{n-1}^{(d)}, \varepsilon_{n}^{(d)}=0$ ($\varepsilon_{i}^{(d)}|\varepsilon_{i+1}^{(d)}$) be the elementary divisors of $T_{S}^{(d)}$.\\ 
\\
{\bf Theorem 5.5} ([58],[59]). {\em For $d\geq 2$, we have}F
$$ e_{d} = \#\{ i \,| \, \varepsilon_{i}^{(d)} \equiv 0 \, {\rm mod} \, \frak{p}^{d} \} - 1.$$
In particular, as an analogue of Cororally 5.3, we have the following formula which generalizes R\'{e}dei's ([73]) for the case $l=d=2$.\\
\\
{\bf Corollary 5.6} ([58],[59]). (1) {\em We have}
$$e_{2} = n - 1 -{\rm rank}_{{\bf F}_{l}}(T_{S} \; {\rm mod} \; l)$$
 {\em where $T_{S} = (T_{S}(i,j))$ is the arithmetic linking matrix defined by} $T_{S}(i,i) = -\sum_{j\neq i}{\rm lk}_{l}(p_{i},p_{j}), T_{S}(i,j)= {\rm lk}_{l}(p_{i},p_{j})$ ($i\neq j$).\\
(2) {\em Suppose that $n =2$. Under the assumption that $e_{d} = 1$ $(d\geq 1)$, we have} 
$$
e_{d+1} = 1 \Longleftrightarrow \left\{
\begin{array}{rl}
\displaystyle{\sum_{r=1}^{d}\sum_{i_{1},\cdots ,i_{r-1} = 1,2}\widehat{\mu}(i_{1}\cdots i_{r-1}21)\pi^{r} \; \equiv \; 0 \; {\rm mod}\; \pi^{d+1}},\\
\displaystyle{\sum_{r=1}^{d}\sum_{i_{1},\cdots ,i_{r-1} = 1,2}\widehat{\mu}(i_{1}\cdots i_{r-1}12)\pi^{r} \; \equiv \; 0 \; {\rm mod}\; \pi^{d+1}}. \\
\end{array}\right.
$$
\\
For the proof, we introduce an arithmetic analogue of the reduced Alexander matrix for the Galois group  $G_{S}(l)$ and translate the whole argument in the case of a link
into our arithmetic situation.\\
\\
{\bf Example 5.7} ([58]). Suppose that $l = 2$ and $S =\{13, 61, 937\}$. By Example 4.7, we have $\mu_{2}(ij) = 0$ ($1\leq i, j \leq 3$), $\mu_{2}(ijk) = 1$ for any permutation $ijk$ of $123$ and $\mu_{2}(ijk) = 0$ otherwise. Furthermore, we see $\mu_{4}(ij) = 0$ ($1\leq i, j \leq 3$). Therefore $T_{S}^{(2)} \equiv  O_{3}$ mod 4 and 
$$  T_{S}^{(3)} \equiv \left( \begin{array}{ccc}
0 & 4  &  4 \\
4 & 0  &  4 \\
4 & 4  &  0 \\
\end{array} \right) \; \sim 
\left( \begin{array}{ccc}
4 & 0  &  0 \\
0 & 4  &  0 \\
0 & 0  &  0 \\
\end{array} \right)\; 
{\rm mod} \; 8.$$
 Hence $e_{2} = 2$, $e_{3} = 0$. By Theorem 5.5, we have $H_{k}(2) = {\bf Z}/4{\bf Z} \oplus {\bf Z}/4{\bf Z}$ for $k = {\bf Q}(\sqrt{13\cdot 61\cdot 937})$.\\

As we mentioned above, our proof of Theorem 5.5 is, as a way of thinking, similar to Iwasawa theory and uses only the classic tools  in the arithmetic of pro-$p$ extensions.  Accordingly, it is a result that should have been obtainable in the 1960's. The reasons for it's going unnoticed seems to be the fact that, at that time, a number field was considered an analogue of a function field and analogues of Milner's invariants were yet to be considered for function fields. (Recently, Y. Terashima and the author ([61]) introduced a multiple generalization of the tame symbol on a Riemann surface using the Massey products in the Deligne cohomology.)  However, once a prime is regarded as a knot not only does our result follows, but a natural generalization along the original line of thought of Gauss follows.
\\
\\
{\bf 6. 3-manifold coverings and number field extensions}
\\

Any oriented, connected closed 3-manifold is realized as a finite
covering of  $S^{3}$ ramified along a link (J. W. Alexander),  just as any number field is a finite extension of ${\bf Q}$ ramified over a finite set of prime numbers. Hence we may expect that there would be conceptual analogies between 3-manifold coverings and number field extensions. First, let us remind a dictionary of the basic analogies between 3-manifolds and number rings ([8],[28],[55],[56],[59],[70],[71]):

$$\begin{array}{ccc}
\mbox{closed 3-manifold}\; M & \longleftrightarrow & \mbox{number ring}\; {\rm Spec}({\cal O}_{k})\cup \{v|\infty \}   \\
 S^{3}  & & {\rm Spec}({\bf Z})\cup \{ \infty \} \\
& & \\
\mbox{knot} \; K  & \longleftrightarrow & \mbox{prime}\; \frak{p} \\
\mbox{link} \, L=K_1\cup\cdots\cup K_n  & \longleftrightarrow & \mbox{primes}\; S = \{ \frak{p}_1,\cdots ,\frak{p}_n\} \\
\mbox{end}& & \mbox{infinite primes}\\
& & \\
\mbox{tubular n.b.d}\; V_{K} & \longleftrightarrow & \frak{p}\, \mbox{adic integer ring} \; {\rm Spec}({\cal O}_{\frak{p}})\\
\mbox{boundary}\; \partial V_{K} & \longleftrightarrow & \frak{p}\, \mbox{adic field} \;  {\rm Spec}(k_{\frak{p}})\\
& & \\
C_{2}(M) \rightarrow C_{1}(M); \; \Sigma \mapsto \partial \Sigma & \longleftrightarrow & k^{\times} \rightarrow I_{k} = \displaystyle{\bigoplus_{\frak{p}}{\bf Z}}; \; a \mapsto (a)   \\
& & \\
H_{1}(M) & \longleftrightarrow & \, H_{k}  \\
& & \\
H_{2}(M) &   \longleftrightarrow & {\cal O}_{k}^{\times} \\
& &\\
\pi_{1}(M) & \longleftrightarrow & \pi_{1}({\rm Spec}({\cal O}_{k}))\\
& &\\
\pi_{1}(M \setminus L)& \longleftrightarrow & \pi_{1}({\rm Spec}({\cal O}_{k})\setminus S) \\
& & \\
\mbox{Hurewicz isomorphism} & \longleftrightarrow & \mbox{unramified classfield theory} \\
 H_{1}(M) \simeq \mbox{Gal}(M^{a}/M)& &  H_{k} \simeq \mbox{Gal}(k^{a}/k) \\
(M^{a}:= \mbox{max. abelian cover of} \, M)& & (k^{a}:= \mbox{Hilbert classfield of} \, k)
\end{array}$$

We note that these analogies are conceptual and we do not claim that
there is a precise one to one correspondence between knots and primes,
3-manifolds and number fields. For example, one has $\pi_{1}({\rm Spec}({\cal O}_{k})) =1$ for any imaginary quadratic field $k$ with class number 1. As an  analogy of the Poincar\'{e} conjecture for number rings, we have the following: \\
{\it Suppose that  $k$ is  a number field whose ring of integers ${\cal O}_k$ 
is ``cohomologically ${\bf Z}$", namely} $ H_{c}^{i}({\rm Spec}({\cal O}_{k}),{\bf Z}) = H_{c}^{i}({\rm Spec}({\bf Z}), {\bf Z}) \; \mbox{for} \, i \geq 0.$ {\em Then ${\cal O}_k$ must be ${\bf Z}$,} \\
where $H_c^{*}$ stands for the compactly-supported \'{e}tale cohomology taking the infinite prime into account ([47,II]). This assertion follows from the Artin-Verdier duality (classfield theory) and Dirichlet's unit theorem. Here we notice that Dirichlet's unit theorem is essentially of analytic nature. According to our analogy, the difficulty of the original Poincar\'{e} conjecture may be coming from that of the corresponding analytic method--gauge theory--in topology.\\

In the following, in light of the above dictionary, we shall discuss analogies of the three classical theories, Hilbert theory, capitulation problem and class towers, in the context of 3-manifold topology. \\
\\
{\bf 6.1.  Hilbert theory.} The Hilbert theory describes the structure of the decomposition of a prime in an extension of number fields. Similarly, we have a topological analogue of the Hilbert theory, which describes the structure of the decomposition of a knot in a 3-manifold covering ([59]).\\

 Let $f : N \rightarrow M$ be a covering of closed 3-manifolds of
finite degree $n$ ramified over a link $L$. Let $K$ be a knot in $M$
 which is a component of $L$ or disjoint from $L$ and let $f^{-1}(K) =
K_{1} \cup \cdots \cup K_{r}$,  a link of $r = r(K)$ components. Let
$V_K$ be a tubular neighborhood of $K$ and let $V_{i}$ be the connected
component of $f^{-1}(V_K)$ containing $K_{i}$.  Take a base point $b
\in \partial V_K$.  We let $f^{-1}(b) = \{ b_{1},\cdots , b_{n} \}$ and consider the monodromy representation $ \psi : G_L = \pi_{1}(M \setminus L, b) \rightarrow {\rm Aut}(f^{-1}(b))$.  Take a finite covering $g : \bar{N} \rightarrow N$  such that the composite $f \circ g : \bar{N} \rightarrow  M$ is a Galois
covering ramified over $L$, and set $G = {\rm Gal}(\bar{N}/M), H = {\rm
Gal}(\bar{N}/N)$. Then $\psi$ factors through $G$ and is identified with the
permutation representation of $G$ on $H\backslash G = \{
 H\sigma_{1},\cdots , H\sigma_{n} \}$. For a knot $\bar{K}$ in $\bar{N}$
 over $K$,  we define the  decomposition group of $\bar{K}$ by
$$ D(\bar{K}) := \{ \sigma \in G \; \vert \; \sigma(\bar{K}) = \bar{K}
\}.$$
Each $ \sigma \in D(\bar{K})$ induces a covering transformation of
$\bar{K}$ over $K$ and we have a homomorphism $D(\bar{K}) \ni \sigma \mapsto \sigma|_{\bar{K}} \in 
{\rm Gal}(\bar{K}/K)$. We define the  inertia group of $\bar{K}$ by

$$I(\bar{K}) := \{ \sigma \in D(\bar{K}) \; \vert \;  \sigma |_{\bar{K}}
= id_{\bar{K}} \}.$$
For another choice of $\bar{K}$ over $K$, $G(\bar{K})$ and $I(\bar{K})$
 are changed up to conjugates. The set of orbits ${\cal O}_{j}$ of
$f^{-1}(b)$ under the action of $D(\bar{K})$ via $\psi$ is identified with
the set of knots $K_{i}$'s:
$$ f^{-1}(b)/D(\bar{K}) \simeq \{ K_{1} ,\cdots, K_{r} \}, \;\; \psi(D(\bar{K}))(b_j) \mapsto q(\sigma_{j}(\bar{K})). $$
In a geometric term, the orbit ${\cal O}_j$ corresponding to $K_j$ is given by the set of
$b_i$ such that $b_{i} \in \partial V_{j}$. We define $f(K_{j})$ by the covering degree of $K_{j}$ over $K$ and the ramification index $e(K_{j})$ of $K_{j}$ by by $\#{\cal O}_{j}/f(K_{j})$. Then we have
the basic identity
$$ n =
\displaystyle{\sum_{j=1}^{r}e(K_{j})f(K_{j})}.$$
  For a Galois covering $f : N \rightarrow M$ with Galois group $G$, the ramification indices and covering degrees are independent of $K_{j}$ and denoted by $e(K)$ and $f(K)$ respectively. \\
\\
{\bf Theorem 6.1.1} ([55],[59],[76]). {\em Let $f : N \rightarrow M$ be a Galois covering. Notations being as above, we have the exact sequence}
$$ 1 \longrightarrow I(K_{j}) \longrightarrow D(K_{j}) \longrightarrow {\rm Gal}(K_{j}/K) \longrightarrow 1,$$
where  $\#I(K_{j}) = e(K)$, $\#D(K_{j}) = e(K)f(K)$, $n = e(K)f(K)r(K)$.\\

Now suppose that $\pi_{1}(M \setminus L)$ is generated normaly by meridians of the components of  $L$. This condition is satisfied when $M$ is a homology 3-sphere (Wirtinger presentation). Then we have the following analogue of the arithmetic progression theorem.\\
\\
{\bf Theorem 6.1.2} ([55],[76]). For any conjugacy class $C$ of ${\rm Gal}(N/M)$, there are infinitely many  knots $K$ (up to ambient isotopy) in $M\setminus L$  so that the conjugacy class of the image of $K$ in ${\rm Gal}(N/M)$ under the map $\pi_1(M\setminus L) \rightarrow {\rm Gal}(N/M)$ is equal to $C$.\\
\\
{\bf 6.2. Capitulation problem.} For a given extension $F/E$ of number fields, the capitulation problem is concerned about the structure of the kernel of the homomorphism $H_E \rightarrow H_F$ asscociated to the extension of ideals. It came originally from the principal ideal theorem, due to Ph. Furtw\"{a}ngler,
  which asserts that any ideal class of a number field becomes principle
in the Hilbert classfield. Since the problem is translated into a purely
group-theoretical problem on the transfer by means of the Artin
reciprocity, we can discuss 3-manifold analogues of
capitulation theorems by means of the Hurewicz isomorphism.\\

Let $N \rightarrow M$ be a finite abelian covering of
path-connected topological manifolds and let $t_{N/M} : H_{1}(M,{\bf Z}) \rightarrow H_{1}(N,{\bf Z})$ be the transfer map. A topological analogue of the capitulation problem is concerned with the abelian group $C_{N/M} := {\rm Ker}(t_{N/M})$.  Let $M^{a}$ and $N^{a}$ be the maximal abelian coverings of $M$ and $N$
 respectively. Then we have the commutative diagram
$$ \begin{array}{cccc}
H_{1}(M,{\bf Z}) & \stackrel{\sim}{\longrightarrow} & \Gamma/\Gamma' &
\;\; \Gamma = {\rm Gal}(N^{a}/M) \\
t_{Y/X}\downarrow & & \downarrow V_{\Gamma/H} \\
H_{1}(N,{\bf Z}) & \stackrel{\sim}{\longrightarrow} & H/H' & \;\;  H =
{\rm Gal}(N^{a}/N)
\end{array} $$
where the horizontal maps are Hurewicz isomorphisms and $V_{\Gamma/H}$
 is the group-transfer. Hence we have
$$ C_{N/M} \simeq {\rm Ker}(V_{\Gamma/H}).$$
Applying a fundamental theorem by H. Suzuki ([81]) on the transfer-kernel, we have the following:\\
\\
{\bf Theorem 6.2.1} ([53]). {\em Let $N \rightarrow M$ be an abelian covering
 of path-connected compact manifolds of degree $d$.
Assume that $H_{1}(M)$ is finite. Then the order of $C_{N/M}$ is divisible by $d$.}\\
\\
In particular, we have a topological analogue of the principal ideal theorem. Theorem 6.2.1 can be generalized to the case where $H_1(M)$ is possibly infinite
([12]).\\

 Let  $f : N
\rightarrow M$ be a finite abelian covering of closed 3-manifolds with Galois group $G$. By using the Poincar\'{e} duality $C_{N/M} \simeq {\rm Ker}(f^{*} : H^2(M) \rightarrow H^2(N))$ and the spectral sequence
$H^{i}(G, H^{j}(N,{\bf Z})) \Longrightarrow
H^{i+j}(M,{\bf Z})$, we can show the following:\\
\\
{\bf Theorem 6.2.2} ([53]). (1) {\em If $N$ is a rational homology 3-sphere, we have an isomorphism}
$$ C_{N/M} \; \simeq \; G/G' \; (\mbox{non-canonical}).$$
(2) {\em Assume that $M$ is a rational homology 3-sphere and
$G$ is cyclic. Then we have the equality}
$$ \#C_{N/M} = [N : M]h(H_{2}(N,{\bf Z})) $$
{\em where $h(H_{2}(N,{\bf Z})) :=
\displaystyle{\prod_{i=0}^{1}\#H^{i}(G,H_{2}(N,{\bf Z}))^{(-1)^{i+1}}}$}
\\
\\
Since $H_2(N)$ is an analogue of a unit group, Theorem 6.2.2, (2) is seen as an analogue of the classical Chevalley-Tate theorem.\\
\\
{\bf 6.3. Class towers.}  The class field tower problem asks whether there is a finite number field $k$ such that infinite tower of Hilbert $l$-classfields over $k$ exists, equivalently, the maximal pro-$l$ quotient $\pi_1({\rm Spec}({\cal O}_k))(l)$ of $\pi_1({\rm Spec}({\cal O}_k))$ is infinite. An affirmative answer was given by Golod and \v{S}afarevi\v{c} ([15]), and their method and its variants are also applicable to 3-manifold groups. In the following, $l$ stands for a prime number.\\

The essential part in Golod and \v{S}afarevi\v{c}'s argument is to show that\\
 {\em  the equality $\dim H^2(G, {\bf F}_l) > \dim H^1(G,{\bf F}_l)/4$ holds for a finite $l$ group $G$.}\\
 It follows from this that $\pi_1({\rm Spec}({\cal O}_k))(l)$ is  an infinte (non-analytic) pro-$l$ group if the inequality 
$$r(k) \geq 2+2\sqrt{s(k)+1} \leqno{(6.3.1)}$$
holds  for the $l$-ranks $r(k)$ and $s(k)$ of $H_k$ and ${\cal O}_k^{\times}/({\cal O}_k^{\times})^l$ respectively (It has been conjectured by Fontaine and Mazur that $\pi_1({\rm Spec}({\cal O}_k))(l)$ is non-analytic whenever it is infinite). To discuss an analogy of the class tower problem over a 3-manifold $M$,
we may assume that $M$ is a rational homology 3-sphere. In
particular, we have $H_{2}(M,{\bf Z}) = 0$ and an analogy of
Golod-\v{S}afarevi\v{c}'s inequality (6.3.1) for $M$ becomes
$\dim_{{\bf F}_{l}}H_{1}(M,{\bf F}_{l}) \geq 4$ for a prime number $l$ in view of the analogy
 between $H_{2}(M,{\bf Z})$ and the unit group of a number field. In fact,
 under this Golod-\v{S}afarevi\v{c}'s condition, it is shown that $\pi_{1}(M)$ has infinitely decreasing lower central
series and enjoys the exponential word growth
property ([75],[85]), in particular, the pro-$l$ completion of $\pi_{1}(M)$ is not $l$-adic analytic. Futhermore, Reznikov has obtained deeper theorems on
the exponential subgroup growth conjecture for 3-manifold groups, which was conjectured by A. Lubotzky for hyperbolic 3-manifolds, by reducing the problem to Artin's primitive root problem.\\
\\
{\bf Theorem 6.3.2} ([70],[72]). (1) {\em Let $M$ be a 3-manifold whose any finite covering is a rational homology 3-sphere and let $M = M_{1} \leftarrow M_{2}
\leftarrow \cdots \leftarrow M_{n} \leftarrow \cdots$ be the tower of
maximal abelian $l$-coverings. Write $r(M_{n})$ for the $l$-rank of
$H_{1}(M_{n},{\bf Z})$ and assume $r(M) \geq 4$. Then $M$ has the infinite
$l$-class tower. More precisely, we have the inequality
$$r(M_{n+1}) \geq \frac{1}{2}(r(M_{n})^{2} - r(M_{n}))$$
 and the pro-$l$ completion of $\pi_{1}(M)$ is an infinite and
non-analytic pro-$l$ group.}
\\
(2) {\em Assume that $M$ hyperbolic or that $M$ is a rational homology 3-sphere whose Casson invariant is greater than the number of representations of
$\pi_{1}(M)$ to $SL_{2}({\bf F}_{5}))$. Then there are positive constants 
$c_1, c_2$ such that for infinitely many $d$, the number of
subgroups of index $d$ in  $\pi_{1}(M)$ is greater than $\exp(c_1
d^{c_2})$.}\\
\\
We note that analogues of Reznikov's theorem for number fields are completely unknown. This is the area remained to be explored in classical algebraic number theory (cf. [3]).
\\
\\
{\bf 7. Alexander-Fox, torsion theory and Iwasawa theory}\\

Based on the analogy between a ${\bf Z}_p$-extension of a number field and the infinite cyclic covering of a knot complement, there are close paralells between Alexander-Fox theory and Iwasawa theory. Furthermore, the Iwasawa main conjecture, which relates Iwasawa polynomials with $p$-adic analytic zeta functions, may be seen as an analogue of the relation between Reidemeister torsions and  analytic torsions.\\

Let $K$ be a knot in a rational homology 3-sphere $M$ and $X = M \setminus K$. Since $H_1(X)=\langle \alpha \rangle \oplus (\mbox{finite group})$, we take the infinite cyclic covering $X_{\infty} \rightarrow X$, ${\rm Gal}(X_{\infty}/X) = \langle \alpha \rangle$, corresponding to $\langle \alpha \rangle = {\bf Z}$. For simplicity, we assume that we can take $\alpha$ to be a meridian of $K$. The homology group $H_1(X_{\infty})$ is then a module over the Laurent polynomial ring $\Lambda := {\bf Z}[{\rm Gal}(X_{\infty}/X)]={\bf Z}[t^{\pm 1}]$ $(\alpha \leftrightarrow t)$ and called the knot module of $K$. The $i$-th elementary ideal $J_i$ of the $\Lambda$-module $H_1(X_{\infty})$ is called the $i$-th Alexander ideal and a generator of its divisorial hull $\cap (J_i)_{\frak{p}}$ ($\frak{p}$ running over prime ideals of $\Lambda$ of height 1) is called the $i$-th Alexander polynomial (cf. [21]). For the case that $M=S^3$, the ideal $J_i$ is described in terms of the Fox free calculus as follows: Given a Wirtinger presentation $G_K = \langle x_1,\dots , x_n \, | \, r_1=\cdots = r_{n-1}=1\rangle$, let $F=\langle x_1,\dots ,x_n \rangle \stackrel{\pi}{\rightarrow} G_K \stackrel{\psi}{\rightarrow} \langle t \rangle$ be the natural map. Then $J_{i-1}$ is given by the ideal generated by $(n-i)$-minors of the Alexander matrix $(\psi\circ \pi (\partial r_i/\partial x_j))$ ($i\geq 1$). In particular, a generator $\Delta_K(t)$ of $J_0$ is called the Alexander polynomial of $K$ and also given as the characteristic polynomial of the action of the meridian $\alpha$ on the knot module $H_1(X_{\infty})$:
$$ \Delta_K(t) = \det(t\cdot id - \alpha \, | \, H_1(X_{\infty})\otimes_{\bf Z}{\bf Q}). \leqno{(7.1)}$$
Denoting by  $M_n$ the Fox completion of the $n$-fold cyclic subcovering of $X_{\infty} \rightarrow X$, we have $H_1(M_n) = H_1(X_{\infty})/(\alpha^n -1)H_1(X_{\infty})$ and hence the knot module $H_1(X_{\infty})$ controls all $H_1(M_n)$ for $n \geq 1$. In particular, the order of $H_1(M_n)$ is equal to $\prod_{\zeta^n=1}|\Delta_K(\zeta)|$ ($\infty$ if $\Delta_K(\zeta)=0$ for a $\zeta$). From this, we obtain the asymptotic formula $\#H_1(M_n) = ab^n$ ($n>>0$; $\log b$ is the Mahler measure of $\Delta_K(t)$) if $\Delta_K(\zeta) \neq 0$ for any $n$-th root $\zeta$ of unity, $n\geq 1$ ([66],[77]). Further, for the $p$-primary part of $H_1(M_{p^{n}})$ where $p$ is a given prime number,  we have $\# H_1(M_{p^n})(p) = p^{\lambda n+\mu p^n + \nu}$ ($n>>0$; $\lambda = {\rm deg}(\Delta_K)$, $\mu, \nu$ being constants) ([21],[27]). For the case $M=S^3$, we see $H_1(M_{p^{n}})(p)=0$ for any $n\geq 0$. \\

On the other hand, let $p$ be a prime number and $k$ a number field of finite degree over ${\bf Q}$, and consider a ${\bf Z}_p$-extension $k_{\infty}$ of $k$ with ${\rm Gal}(k_{\infty}/k) = \langle \gamma \rangle = {\bf Z}_p$. For simplicity, we assume that there exists only one prime $\frak{p}$ of $k$ lying over $p$ and that $\frak{p}$ is totally ramified in $k_{\infty}/k$. Denote by $k_n$ the cyclic subextension of $k_{\infty}/k$ of degree $p^n$ over $k$ and by $H_n$ the $p$-Sylow subgroup of the ideal class group of $k_n$ for $n\geq 0$. The Iwasawa module of $(p,k)$ is then defined by $H_{\infty} := \varprojlim_{n} H_n$ (projective limit being taken with respect to the norm maps), which is regarded as a module over the Iwasawa algebra $\hat{\Lambda} := {\bf Z}_p[[{\rm Gal}(k_{\infty}/k)]]={\bf Z}_p[[T]]$ $(\gamma \leftrightarrow 1+T)$. The $i$-th elementary ideal of the $\hat{\Lambda}$-module $H_{\infty}$ is called the $i$-th Iwasawa ideal and a polynomial generator of its divisorial hull is called the $i$-th Iwasawa polynomial. As in the case of a knot, the $i$-th Iwasawa ideal may be described in terms of the pro-$p$ Fox free calculus by using a presentation of the pro-$p$ Galois group ${\rm Gal}(\tilde{k}_{\infty}/k)$ of the maximal unramified pro-$p$ extension $\tilde{k}_{\infty}$ of $k_{\infty}$ over $k$. In particular, as an analogue of $(7.1)$, the $0$-th Iwasawa polynomial is given by the characteristic polynomial of the action of monodromy $\gamma$ on the Iwasawa module $H_{\infty}$ (up to a power $p^{\mu}$ of $p$):
$$ I_p(T) = \det(T\cdot id -(\gamma -1)\, | \, H_{\infty}\otimes_{{\bf Z}_p}{\bf Q}_p). \leqno{(7.2)}$$
In view of the analogy with the function field of an algebraic curve over a finite field, $\gamma$ is regarded as an analogue of the Frobenius automorphism of the unramified extension arising from the constant field extension. We note, however, that it is more natural to regard $\gamma$ 
 as an analogue of the meridian $\alpha$, because $\gamma$ is a generator of the inertia group of the ramified extension $k_{\infty}/k$ over $p$. As in the case of a knot, the Iwasawa module $H_{\infty}$ controls all $H_n$ for $n\geq 0$, namely we have $H_n = H_{\infty}/(\gamma^{p^n}-1)H_{\infty}$, and in particular the Iwasawa formula $\# H_n = p^{\lambda n + \mu p^n + \nu}$ ($n>>0$; $\lambda$ = {\rm deg}($I_p(T)$), $\mu, \nu$ are constants). For $k={\bf Q}$, we have $H_n=0$ for any $n \geq 0$. In particular, the Iwasawa $\lambda$-invariant is seen as an analog of twice of the genus of a knot in case that a knot is fibred.\\

In fact, we can construct the Alexander (resp. Iwasawa) module as the $\psi$-diffrential module for any group (resp. pro-$p$ group) $G$ and a homomorphism $\psi : G \rightarrow {\bf Z}$ (resp. ${\bf Z}_p$) ([20],[51]), and the analogies discussed above are consequences obtained in this general group-theoretic framework. Furthermore, the both theories are generalized to those on the twisted Alexander module (resp. twisted Iwasawa module, called usually Selmer group) associated to a complex representation of a knot group (resp. $p$-adic representation of a Galois group) ([16],[31]). In both theories, however, the deeper aspect is the connection of these algebraic invariants with ($p$-adic) analytic zeta functions.\\

For a knot $K \subset S^3$, set $X := S^3 \setminus V_K^o$ and $G_K := \pi_1(X)$. Given a finite dimensional representation $\rho : G_K \rightarrow GL(V)$ over a field $F$, let $C_*(X,\rho)$ be the chain complex with coefficient $\rho$ and $H_*(X,\rho)$ the homology groups. An important role in the theory of torsions is played by the Euler isomorphism $Eu : \det C_{*}(X,\rho) \stackrel{\sim}{\rightarrow} \det H_{*}(X,\rho)$ between the determinant modules ([32]). If $C_{*}(X,\rho)$ is acyclic, we have $\det C_{*}(X,\rho) \simeq F^{\times}$. Choosing a basis $b = \{ b_{i} \}$ of $C_{*}(X,\rho)$, the Reidemeister torsion $\tau(X,\rho) := Eu(\det(b))$ is defined up to $\pm \det (\rho(G_{K}))$ and gives a topological invariant of $(X,\rho)$. For the case that $F = {\bf C}(t)$ and  $\rho = \psi : G_{K} \rightarrow \langle t \rangle \subset {\bf C}(t)^{\times}$ is the universal family of 1-dimensional representations, $C_{*}(X,\psi)$ is acyclic and $\tau(X,\rho)$ is essentially given by the Alexander polynomial $\Delta_K(t)$ of $K$. In the following, we assume that $F = {\bf R}$ and $\rho$ is an orthogonal representation. Since $X$ is a compact manifold with boundary,  we define the Reidemeister torsion $\tau = \tau(X,\partial X;\rho)$ by using the relative chain complex $C_{*}(X,\partial X;\rho)$ (well defined up to $\pm 1$). Let $\rho^{*}$ be the adjoint representation of $\rho$ and $E_{\rho^{*}}$   the associated flat vector bundle over $X$. We may assume that the Riemannian metric of $X$ is given as a direct product near $\partial X$. Using the metrics on $X$ and $E_{\rho^{*}}$, we regard $\Omega^{i}(X,\rho^{*})$ of $E_{\rho^{*}}$-valued $i$-forms as a Hilbert space. Set $\Omega^{i}(X,\partial X;\rho^{*}) = \{ \omega \in \Omega^{i}(X,\rho^{*}) \,| \, \omega|_{\partial X} = \delta\omega|_{\partial X} = 0 \}$ ($\delta$ being the adjoint operator of the differential $d$) and we define the spectral zeta function of the Laplace operator $\Delta^{i}$ restricted on  $\Omega^{i}(X,\partial X;\rho^{*})$ by $Z(\rho^{*},s) := \sum_{i=0}^{3}(-1)^{i+1}i\zeta_{i}(\rho^{*}, s)$, $\zeta_{i}(\rho^{*},s) := \sum \lambda_{i}^{-s}$, where $\lambda_{i}$ runs over positive eigenvalues of  $\Delta^{i}|_{\Omega^{i}(X,\partial X;\rho^{*})}$. It is shown that $\zeta_{i}(\rho^{*},s)$ is continued meromorphically to the ${\bf C}$-plane and that it is holomorphic at $s=0$. Via the de Rham isomorphism with the space of harmonic forms, $H^{i}(X,\partial X;\rho)$ is equipped with a norm $|\,\cdot\,|$. Take a basis $\xi \in \det H^{i}(X,\partial X;\rho)$ such that $|\xi| = 1$ and we define the zeta element $z(X,\rho)$ by $\tau\otimes \xi$. Then J.Cheeger-W.M\"{u}ller-W.L\"{u}ck's theorem ([39]) is stated as follows:  $z(X,\rho)$ is mapped to $\pm \exp(\frac{1}{2}Z'(\rho^{*},0))$ under the isomorphism $\det H_{*}(X,\partial X;\rho) \otimes \det H^{*}(X,\partial X;\rho^{*}) \stackrel{\sim}{\rightarrow} {\bf R}$.\\

On the other hand, the Iwasawa polynomial is also related with a ($p$-adic) analytic zeta function. This connection is called the Iwasawa main conjecture. For a prime number $p > 2$, let $X = {\rm Spec}({\bf Z})\setminus \{(p)\}$ and $G_{\{ (p)\}} = \pi_{1}(X)$. Let $\psi : G_{\{(p)\}} \rightarrow {\rm Gal}({\bf Q}(\mu_{p^{\infty}})/{\bf Q}) \subset \tilde{\Lambda}^{\times}$ be the universal family of 1-dimensional $p$-adic representations of $G_{\{(p)\}}$ where $ \tilde{\Lambda} := {\bf Z}_{p}[[{\rm Gal}({\bf Q}(\mu_{p^{\infty}})/{\bf Q})]]$ and  $\mu_{p^{\infty}}$ stands for the group of all $p^n$-th roots of unity for $n\geq 1$. As in the case of the Alexander polynomial, the Iwasawa polynomial is regarded as a part of the determinant module $\det H^{*}(X,\psi)\otimes_{\tilde{\Lambda}}Q(\tilde{\Lambda})$ ($Q(\tilde{\Lambda})$ being the total quotient ring of $\tilde{\Lambda}$). K. Kato ([30]) gives a torsion type formulation of the generalized Iwasawa main conjecture as follows. Let $M$ be a motif $H^{m}(Y)(r)$ over ${\bf Q}$ ($Y$ being a smooth proper scheme over ${\bf Q}$). Assume that the action of ${\rm Gal}(\bar{\bf Q}/{\bf Q})$ on $M_{p} = H^{m}(Y\otimes \bar{\bf Q}, {\bf Q}_{p})(r)$ is unramified outside $\{p,\infty\}$ and let $\rho$ be the associated $p$-adic representation  of $G_{\{ (p)\}}$, and consider $\Delta (X,\rho) := \det {\rm R}\Gamma(X,M_{p})^{-1}\otimes_{{\bf Q}_{p}}\det {\rm R}\Gamma(X\otimes {\bf R},M_{p}(-1))^{-1}$. After the model of the theory of Euler systems, K. Kato ([ibid]) conjectured that there is a good basis $z(X,\rho)$ of $\Delta(X,\rho)$ which is related with the zeta function  $L(\rho^{*}(1),s) := \prod_{l\neq p} \det(I-{\rm Fr}_{l}^{-1} l^{-s} | M_{p}^{*}(1))^{-1}$ attached to $M_{p}^{*}(1) = {\rm Hom}(M_{p},{\bf Q}_{p})(1)$. This zeta element $z(X,\rho)$, which is which living in $\Delta (X,\rho)$, is seen as a $p$-adic analogue of the Reidemeister torsion. Let $\Sigma := H^{m}(Y({\bf C}),{\bf Q}(2\pi i)^{r-1})^{+}$ ($+$ means the invariant part under the action of complex conjugation) and  $\Omega := {\rm Im}(H^{m}(Y,\Omega_{Y/{\bf Q}}^{\geq r}) \rightarrow H^{m}_{\rm dR}(Y/{\bf Q}))$, and assume that  the ($p$-adic) period maps $a : \Omega \otimes {\bf R} \rightarrow \Sigma \otimes {\bf R}$ and  $b : H^{1}(X,M_{p}) \rightarrow \Omega\otimes {\bf Q}_{p}$ are isomorphisms and that $H^{i}(X,M_{p}) = 0$ ($i \neq 1$) (strictly critical case [30, $\S$3.2$\sim$3.3]). Then we have  $\Delta(X,\rho) = \det H^{1}(X,M_{p})\otimes_{\bf Q} \det \Sigma^{-1}$ and  isomorphisms $f : \det(\Omega)\otimes_{\bf Q} \det(\Sigma)^{-1}\otimes {\bf R} \stackrel{\sim}{\rightarrow} {\bf R}$ and $g : \det(\Omega)\otimes_{\bf Q}\det (\Sigma)^{-1}\otimes {\bf Q}_{p} \stackrel{\sim}{\rightarrow} \Delta(X, \rho)$. Then Kato's conjecture is stated as follows: $z(X,\rho)$ is mapped to an element of  $\det(\Omega)\otimes \det (\Sigma)^{-1}$ under $g^{-1}$ and further mapped to $\pm L(\rho^{*}(1),0)$ under $f$. \\

As we observed above, there is a close parallel between the frameworks of torsion theory and Iwasawa theory, and hence both theories could provide the ideas each other.  In recent years, M. Kurihara ([36]) studied the higher elementary ideals of the Iwasawa module and showed a refined form of the Iwasawa main conjecture. After the works of Milnor ([49]) and D. Fried ([10]), K. Sugiyama ([78],[79]) also showed a geometric analogue of the Iwasawa main conjecture for a hyperbolic 3-manifold having an infinite cyclic cover. It should also be mentioned that torsion theory is generalized to an equivariant theory (where a manifold has a group action) using the $K$-groups of group von Neumann algebras ($L^2$-torsion theory), while equivariant versions of the Bloch-Kato conjecture and Iwasawa theory are also studied using the $K$-groups of completed group algebras (cf. [4],[40]).\\
 
{\bf 8. Moduli and deformations of representations of knot and prime groups} \\

From the viewpoint of the analogy between the knot group $G_{K} = \pi_{1}(S^{3}\setminus K)$ and the prime group  $G_{\{(p)\}} = \pi_{1}({\rm Spec}({\bf Z}) \setminus \{(p)\})$, one may expect that there would also be analogies between the moduli(deformation) spaces of representations of $G_K$ and $G_{\{ (p)\}}$. In the following, we discuss some analogies on these moduli(deformation) spaces and invariants defined on them. \\

For a knot $K \subset S^3$ and a tubular neighborhood $V_{K}$ of $K$, we choose  a meridian $\alpha$ and a longitude $\beta$ on $\partial V_K$ and set $I_{K} := \langle \alpha \rangle \subset D_{K} := \pi_{1}(\partial V_{K})$. On the other hand, for a prime number $p$, we set $I_{\{(p)\}} := \mbox{the inertia group over $p$} \; \subset D_{\{ (p) \}} := \pi_{1}({\rm Spec}({\bf Q}_{p}))$. According to (1.6), the restriction of a representation of $G_{K}$ (resp. $G_{\{(p)\}}$) to $D_K$ (resp. $D_{\{ (p)\}}$) gives an important information as a ``boundary condition".\\

First, consider the set of 1-dimensional complex representations of the knot group $G_{K}$: ${\cal X}_{K}^{1} = {\rm Hom}(G_{K},{\bf C}^{\times})$. Since $G_{K}$ is generated by the conjugates of $I_{K}$, $\rho \in {\cal X}_{K}^{1}$ is determined by the image $\rho(\alpha)$ of a meridian $\alpha$. Thus ${\cal X}_{K}^{1}$ is identified with ${\bf C}^{\times}$. Therefore, the $i$-the Alexander ideals $J_{K,i}$ are regarded as a ascending chain of ideals of the coordinate ring ${\bf C}[{\cal X}_{K}^{1}] = {\bf C}[t^{\pm 1}]$ of ${\cal X}_K^1$ and they define the Alexander stratification  ${\cal X}_{K}(i):= \{ \rho \in {\cal X}_{K}^{1} \, | \, f(\rho(\alpha)) = 0, \; \forall f \in J_{K,i-1} \}$($i\geq 1$) in ${\cal X}_K^1$. Then E. Hironaka ([22]) characterized ${\cal X}_{K}(i)$'s as the jumping loci of the cohomology $H^{1}(G_{K}, \rho)$ as follows. Set ${\cal Y}_{K}(i) := \{ \rho \in {\cal X}_{K}^{1} \, | \, \dim_{\bf C}H^{1}(G_{K}, \rho) \geq  i \}$. Then we have
$$ {\cal Y}_{K}(i) = {\cal X}_{K}(i) \, (i\neq 1), \; {\cal Y}_{K}(1) = {\cal X}_{K}(1) \cup \{ {\bf 1} \} \; ({\bf 1}\; \mbox{\em is the trivial repr.})
\leqno{(8.1)}$$
In particular, We have $\Delta_K(\rho(\alpha))=0 \Leftrightarrow {\rm dim}_{\bf C}H^{1}(G_{K},\rho) > 0$ ([38]).\\

Next, consider 2-dimensional representations of $G_K$. Since $H^{2}(G_{K}, {\bf Z}/2{\bf Z}) = H^{2}(S^{3}\setminus K,{\bf Z}/2{\bf Z}) = 0$, a representation $G_{K} \rightarrow PGL_{2}({\bf C}) = PSL_{2}({\bf C})$ can always be lifted to a representation  to $SL_{2}({\bf C})$. Hence we may consider only representations of $G_K$ to $SL_2({\bf C})$ if the center of $G_{K}$ is trivial (this is the case except torus knots). In the following, we suppose that $S^{3} \setminus K$ is a complete hyperbolic 3-manifold with finite volume (This is satisfied unless $K$ is a torus knot or satellite knot). Let  ${\cal X}_{K}^{2}$ be the character variety of $SL_{2}({\bf C})$-representations of $G_K$ ([6]): ${\cal X}_{K}^{2} = {\rm Hom}(G_{K},SL_{2}({\bf C}))//SL_{2}({\bf C}).$  For $[\rho] \in {\cal X}_{K}^{2}$, as $D_{K}$ is abelian, $\rho|_{D_{K}}$ is equivalent to an upper triangular representation:
$$ \rho|_{D_K} \simeq \left(
\begin{array}{cc}
\chi_{\rho} & * \\
0 & \chi_{\rho}^{-1}\\
\end{array}
\right).
$$
Let $\rho_{o} : G_{K} \rightarrow SL_{2}({\bf C})$ be a lift of the holonomy representation attached to the complete hyperbolic structure on $S^{3}\setminus K$. Then W. Thurston showed the following ([64],[82]):\\
\\
(8.2)$\;$ {\em The map} $\Psi_K\, : \,{\cal X}_{K}^{2} \rightarrow {\bf C}$ {\em defined by}  $\Psi_K([\rho]):= {\rm Tr}(\rho(\alpha))$ {\em is bianalytic in}
{\em a neighborhood  $V$ of} $\rho_{o}$. \\
\\
We set $\chi_{\rho}(\alpha) = \exp(x_K(\rho)), \chi_{\rho}(\beta) = \exp(y_K(\rho))$. Then $x_K, y_K$ are lifted to holomorphic functions, denoted by the same $x_K, y_K$ respectively, on the deformation space $U$ of hyperbolic structures on $S^3 \setminus K$, which is a double covering of $V$ branched at the complete hyperbolic structure $z_{o}$ lying over $\rho_o$ so that $x_K$ gives a local coordinate of $U$ around $z_{o}$. Let $q_K$ be the multiplicative period of $\partial V_K$: $\partial V_K = {\bf C}^{\times}/q_K^{\bf Z}$. Then Neumann-Zagier showed the following relation ([64]):
$$ \frac{dy_K}{dx_K}|_{x_K = z_{o}} = \frac{1}{2\pi\sqrt{-1}} \log(q_K).
 \leqno{(8.3)} $$
We also note that the integral $\int_0^{z} y_K dx_K$ gives (essentially) the $SL_2({\bf C})$ Chern-Simons invariant for the holonomy representation associated to $z$ ([62]).\\

On the other hand, for a prime number $p$, we consider infinitesimal deformations of a given mod $p$ representation $\bar{\rho} : G_{\{ (p) \}} \rightarrow GL_{n}({\bf F}_{p})$. Namely we call a pair $(R,\rho)$ a deformation of $\bar{\rho}$ if $R$ is a complete noetherian local ring with residue field $R/\frak{m}_{R} = {\bf F}_p$ ($\frak{m}_R$ being the maximal ideal of $R$) and $\rho$ is a continuous representaion $G_{\{(p)\}} \rightarrow GL_{n}(R)$ such that $\rho$ mod $\frak{m}_{R} = \bar{\rho}$. For simplicity, we assume in the following that $\bar{\rho}$ is absolutely irreducible. Due to Mazur's fundamental theorem in Galois deformation theory ([45]), there are the universal deformation ring $R_{n}(\bar{\rho})$ and the universal deformation $\rho_{n} : G_{\{ (p) \}} \rightarrow GL_{n}(R_{n}(\bar{\rho}))$ of $\bar{\rho}$ such that any deformation $(R,\rho)$ of $\bar{\rho}$  is given by $\varphi \circ \rho_{n} = P\circ \rho \circ P^{-1}$ ($ P \in I_{n} + M_{n}(\frak{m}_{R})$) for a ${\bf Z}_{p}$-algebra homomorphism $\varphi : R_{n}(\bar{\rho}) \rightarrow R$.  We then  define the universal deformation space of  $\bar{\rho}$ by the rigid analytic space ${\cal X}_{p}^{n}(\bar{\rho}) := {\rm Hom}_{c}(R_{n}(\bar{\rho}),{\bf C}_{p})$ where ${\bf C}_{p}$ be the $p$-adic completion  of an algebraic closure of ${\bf Q}_{p}$. We write $\rho_{\varphi}$ for the representation  $\varphi \circ \rho_{n}$ over ${\bf C}_p$ corresponding to $\varphi \in {\cal X}_{p}^{n}(\bar{\rho})$. In the following, we assume $p>2$ for simplicity.\\

First, consider the deformation space ${\cal X}_{p}^{1}(\bar{\rho})$ of 1-dimensional representations. Let ${\bf Q}_{\infty}/{\bf Q}$ be the ${\bf Z}_{p}$-extension with Galois group $\Gamma = {\rm Gal}({\bf Q}_{\infty}/{\bf Q}) = \langle \gamma \rangle$. Then the universal deformation ring $R_{1}(\bar{\rho})$ is given by the Iwasawa algebra $\hat{\Lambda}:= {\bf Z}_{p}[[\Gamma]] = {\bf Z}_{p}[[T]]$ ($\gamma \leftrightarrow 1+T$) and the universal deformation $\rho_{1} : G_{\{ (p) \}} \rightarrow R(\bar{\rho})^{\times}$ is given by $\rho_{1}(g) = \tilde{\rho}(g)g_{p}$ where $\tilde{\rho}$ is the Teichm\"{u}ller lift of  $\bar{\rho}$ and  $g_{p}$  is the image of $g$ under the natural map $G_{\{ (p) \}} \rightarrow \Gamma$.  Hence ${\cal X}_{p}^{1}(\bar{\rho})$ is identified with the $p$-adic unit disk ${\cal D}_{p}^{1} := \{ x \in {\bf C}_{p} \, | \, |x|_{p} < 1\}$ via the correspondence $\varphi \mapsto \varphi(T) = \varphi(\gamma)-1$. Then we have an analogy of Hironaka's theorem (8.1) as follows. Let $\mu_{p^{n}}$ be the group of $p^{n}$-th roots of unity and $\mu_{p^{\infty}} := \cup_{n\geq 1}\mu_{p^{n}}$. For the Teichm\"{u}ller character $\omega : G_{\{(p)\}} \rightarrow {\bf Z}_p^{\times}$ defined by the action of $G_{\{ (p)\}}$ on $\mu_{p}$, set ${\cal X}_{p}^{1, (j)} := {\cal X}_{p}^{1}(\omega^{j}\, {\rm mod} \,p)$. The Galois group ${\rm Gal}({\bf Q}(\mu_{p^{\infty}})/{\bf Q}) \simeq \Gamma \times {\rm Gal}({\bf Q}(\mu_{p})/{\bf Q})$ acts on $H_{\infty}:= H_{1}({\rm Spec}({\bf Z}[\mu_{p^{\infty}},p^{-1}]),{\bf Z}_{p})$. Let   $H_{\infty}^{(j)}$ be the maximal $\hat{\Lambda}$-submodule of  $H_{\infty}$ on which $\delta \in {\rm Gal}({\bf Q}(\mu_{p})/{\bf Q})$ acts as multiplication by $\omega^{j}(\delta)$. Let $J_{p,i}^{(j)}$ be the $i$-th elementary ideal of $H_{\infty}^{(j)}$ and set ${\cal X}_{p}^{(j)}(i) := \{ \varphi \in {\cal X}_{p}^{1, (j)} | f(\varphi(\gamma)-1) = 0, \; \forall f \in J_{p,i-1}^{(j)} \}$($i\geq 1$), ${\cal Y}_{p}^{(j)}(i) := \{ \varphi \in {\cal X}_{p}^{1, (j)} | \dim_{{\bf C}_{p}} H^{1}(G_{\{(p)\}},\rho_{\varphi}) \geq i \}$. Then for an even integer $j$ ($0\leq j \leq p-2$), we have
$$ {\cal Y}_{p}^{(j)}(i) = {\cal X}_{p}^{(j)}(i)\; (i \neq 1),\;\; 
{\cal Y}_{p}^{(j)}(1) = \left\{
\begin{array}{l}
{\cal X}_{p}^{(j)}(1) \;\; (j\neq 0), \\
 {\cal X}_{p}^{(0)}(1) \cup \{ {\bf 1} \}\;\; (j=0).
\end{array}\right.
\leqno{(8.4)}
$$
 where ${\bf 1}$ is the ${\bf Z}_{p}$-algebra homomorphism defined by ${\bf 1}(\gamma)=1.$ In particular, for the Iwawasa polynomial $I_{p}^{(j)}(T) = \det (T\cdot id - (\gamma-1) | H_{\infty}^{(j)}\otimes_{{\bf Z}_{p}}{\bf Q}_{p})$ attached to $H_{\infty}^{(j)}$, we have $I_{p}^{(j)}(\varphi(\gamma)-1) = 0 \; \Leftrightarrow \; \dim_{{\bf C}_{p}}H^{1}(G_{\{(p)\}},\rho_{\varphi}) > 0$ $(j\neq 0)$ ([57]).\\

Next, let us see an analogy (8.2) for Galois deformation, which was originally pointed out by K. Fujiwara ([11]). In the following, we assume that $\bar{\rho}(\mbox{complex conjugation})$ is not scaler. In order to obtain an analogy of the case of knots, it is natural to put the condition that the restriction of a Galois representation to $D_{\{(p)\}}$  is equivalent to an upper-triangular representation. This ``boundary condition" is called the (nearly) $p$-ordinary condition ([45]). To be precise, a continuous representation $\rho : G_{\{ (p) \}} \rightarrow GL_{2}(R)$ is called  $p$-ordinary if  the $\rho|_{D_{\{ (p) \}}}$ is equivalent to an upper triangular representation:
$$ \rho|_{D_{\{ (p) \}}} \simeq \left(
\begin{array}{cc}
\chi_{\rho,1} & * \\
0 & \chi_{\rho,2}\\
\end{array}
\right), \;\;\;\;\; \chi_{\rho,2}|_{I_{\{ (p) \} }} = {\bf 1}.
$$
It is shown that for a $p$-ordinary residual representation $\bar{\rho}$ there is the universal $p$-ordinary deformation $\rho_{2}^{ord} : G_{\{ (p) \}} \rightarrow GL_{2}(R_{2}^{ord}(\bar{\rho}))$ ([45]). We define the universal $p$-ordinary deformation space by ${\cal X}_{p}^{2,ord}(\bar{\rho}) := {\rm Hom}_{c}(R_{2}^{ord}(\bar{\rho}),{\bf C}_{p})$.  Now we suppose that $\bar{\rho}$ is a residual representation attached to a $p$-ordinary Hecke eigen cuspform on $\Gamma_0(p^{m})$ $(m\geq 1)$. Then it was conjectured by Mazur and proved by A. Wiles etc  under some conditions that $R_{2}^{ord}(\bar{\rho})$ coincides with Hida's  $p$-ordinary Hecke algebra, and we assume it in the following. Then, by H. Hida's theorem ([18],[19]), the natural morphism ${\rm Spec}(R_{2}^{ord}(\bar{\rho})) \rightarrow {\rm Spec}(\hat{\Lambda})$ is finite flat and \'{e}tale at arithmetic points corresponding to modular forms with integral weight. Hence we have the following analogue of (8.2) ([59],[60]):  Let $\rho_{\varphi_o}^{ord} := \varphi_o \circ \rho_{2}^{ord}$ be the $p$-adic representation of $G_{\{(p)\}}$ attached to a $p$-ordinary Hecke eigen cusp form  where $\bar{\rho} = \rho_{\varphi_{o}}^{ord}$ mod $p$. We take $\tau \in I_{\{(p)\}}$, as an analogue of a meridian,  so that $\tau_p = \gamma$. \\
\\
(8.5)$\;$ {\em The map} $\Psi_p \,:\, {\cal X}_{p}^{2,ord}(\bar{\rho}) \rightarrow {\bf C}_{p}$ {\em defined by} $\Psi_p(\varphi):= {\rm Tr}(\rho_{\varphi}^{ord}(\tau))$ {\em is bianalytic in a neighborhood of} $\varphi_o$. \\ 
\\
In order to see an analogy of (8.3), we assume that $\bar{\rho}$ is the mod $p$ representation of the $p$-adic representation $\rho_E$ of $G_{\{ (p)\}}$ arising from the $p$-adic Tate module of a $p$-ordinary, modular elliptic curve $E$ over ${\bf Q}$ which has good reduction outside $p$ and split multiplicative reduction at $p$. We let
$$ \rho_2^{ord}|_{D_{\{ (p)\}}} \simeq \left(
\begin{array}{cc}
\chi_1^{ord} & *\\
0 & \chi_2^{ord}
\end{array}
\right), \;\;\;\;\; \chi_2|_{I_{\{ (p) \}}} = {\bf 1}
$$
and set $x_p = \omega_1^{-1}(\tau) \chi_1^{ord}(\tau)$ and $y_p = \chi_2^{ord}(\sigma)$ where $\omega_1$ is the Teichm\"{u}ller lift of $\det \bar{\rho}$ and $\sigma$ is a Frobenius automorphism over $p$. Then $x_p$ is regarded as a local coordinate of ${\cal X}_p^{2,ord}(\bar{\rho})$ around $\varphi_{E}$ where $\varphi_{E} \circ \rho_2^{ord} = \rho_{E}$. Let $q_E \in p{\bf Z}_p$ be the Tate period of $E$: $E(\overline{{\bf Q}}_p) = \overline{{\bf Q}}_p^{\times}/q_E^{\bf Z}$. Then Greenberg-Stevens showed the following relation ([17]):
$$ \frac{d y_p}{d x_p}|_{x_p = \varphi_E} = -\frac{1}{2} \frac{ \log_p(q_E)}{{\rm ord}_p(q_E)}. \leqno{(8.6)}$$
Here $\log_p$ and ${\rm ord}_p$ stand for the $p$-adic logarithm and the normalized valuation on $\overline{\bf Q}_p$, respectively. The quantity $\log_p(q_E)/{\rm ord}_p(q_E)$ is called the ${\cal L}$-invariant of the elliptic curve $E$.
\\
\\
In view of (8.2), (8.3) and (8.5), (8.6), we find analogies between $(x_K, y_K)$ on ${\cal X}_K^2$ and $(x_p, y_p)$ on ${\cal X}_p^{2, ord}(\bar{\rho})$. It would be interesting to study further analogies between 3-dimensional hyperbolic geometry and the arithmetic of modular Galois deformations. For example, it is interestiong to pursue the analogies between the Chern-Simons functional on the deformation space of hyperbolic structures and the $2$-variable $p$-adic $L$-function on the deformation space of ordinary modular Galois representations (For partial results, see [60], [62]). Furthermore, it would also be an interesting problem to extend the analogy between (8.1) and (8.3) for higher dimensional representations (cf. [46]).\\

Finally, we mention some related works. M. Kapranov ([28]) pointed out the analogies between the structures of topological quantum field theory (M. Atiyah) and Langlands correspondence, where Witten invariants defined as Feynman integrals correspond to automorphic $L$-functions defined as integrals on adele groups. This work is closely related to our analogies in spirit.  On the other hand, M. Kontsevich ([35]) reconstructed the Witten invariants of knots as the universal quantum invariants which turned out to be related, via the Drinfeld associator,  with the arithmetic objects such as multiple zeta values and the absolute Galois group ${\rm Gal}(\overline{\bf Q}/{\bf Q})$ ([2],[9],[24]). It would be an interesting problem  to find a connection between our analogies and the arithmetic aspect appearing in the study of the Kontsevich invariants.

\vspace{.5cm}

{\footnotesize
\begin{flushleft}
{\bf References}\\
{[1]} D. Anick, Inert sets and the Lie algebra associated to a group, 
 J. Algebra, {\bf 111}, (1987), 154--165.\\
{[2]} D. Bar-Natan, On associators and the Grothendieck-Teichm\"{u}ller
   group, I. Selecta Math. (N.S.) {\bf 4}, (1998), 183-212. \\
{[3]} N. Boston, Galois $p$-groups unramified at $p$---a survey, in: Primes and Knots, 31--40, Contemp. Math., {\bf 416}, Amer. Math. Soc., Providence, RI, 2006.\\
{[4]} J. Coates, T. Fukaya, K. Kato, R. Sujatha, O. Venjakob, The ${\rm GL}_2$ main conjecture for elliptic curves without complex multiplication, Publ. Math. Inst. Hautes Etudes Sci. No. {\bf 101}, (2005), 163--208.\\
{[5]} R. Crowell, R. Fox, Introduction to knot theory, Springer GTM, {\bf 57},  New York-Heidelberg, 1977. \\
{[6]} M. Culler, P. Shalen, Varieties of group representations and
   splittings of $3$-manifolds, Ann. of Math. (2) {\bf 117}, (1983), 109-146.\\
{[7]} C. Deninger, Some analogies between number theory and dynamical systems on foliated spaces, In: Proceedings of the International Congress of Mathematicians, Vol. I (Berlin, 1998). Doc. Math. 1998, Extra Vol. I, 163--186. \\
{[8]} --------, A note on arithmetic topology and dynamical systems, in: Algebraic number theory and algebraic geometry, Contemp. Math., {\bf 300}, AMS, (2002), 99--114.\\
{[9]} V.G. Drinfeld, On quasitriangular quasi-Hopf algebras and on
   a group that is closely connected with ${\rm Gal}(\overline{\bf Q}/{\bf Q})$, Leningrad Math. J. {\bf 2}, (1991), no. 4, 829-860.\\
{[10]} D. Fried, Analytic torsion and closed geodesics on hyperbolic manifolds, Invent. Math. {\bf 84}, (1986), no. 3, 523--540. \\
{[11]} K. Fujiwara, $p$-adic gauge theory in number theory, JAMI conference ``Primes and Knots", 2003, March.\\
{[12]} K. Fuluta, On capitulation theorems for infinite groups, in: Primes and Knots, Contemp. Math., {\bf 416}, Amer. Math. Soc., (2006), 41--47.\\
{[13]} C.F. Gauss, Disquisitiones arithmeticae, Yale Univ., (1966).\\
{[14]}--------, Zur mathematischen Theorie der electrodynamischen Wirkungen, Werke {\bf V}, (1833), p.605.\\
{[15]}E. Golod, I. \v{S}afarevi\v{c}, On the class field tower,
(Russian) Izv. Akad. Nauk
SSSR Ser. Mat. {\bf 28}, (1964) 261--272.\\
{[16]} R. Greenberg, Iwasawa theory for $p$-adic representations, in: Algebraic number theory, Adv. Stud. Pure Math., {\bf 17}, Academic Press, Boston, MA, (1989), 97--137. \\
{[17]}  R. Greenberg, G. Stevens, $p$-adic $L$-functions and $p$-adic periods of modular forms, Invent. Math. {\bf 111}, (1993), no. 2, 407--447.\\
{[18]} H. Hida, Iwasawa modules attached to congruences of cusp forms, Ann. Sci. \'{E}cole Norm. Sup. (4) {\bf 19}, (1986), 231--273.\\
{[19]} --------, Galois representations into ${\rm GL}_2({\bf Z}_{p}[[X]])$ attached to ordinary cusp forms, Invent. Math. {\bf 85}, (1986), 545--613.\\
{[20]} J. Hillman, Algebraic invariants of links, Series on Knots and
   Everything, {\bf 32}, World Sci. Publ. Co., 2002.\\
{[21]}J. Hillman, D. Matei, M. Morishita, Pro-$p$ link groups and $p$-homology groups, in: Primes and Knots, Contemporary Math., {\bf 416},  Amer. Math. Soc., (2006), 121-136.\\
{[22]} E. Hironaka,  Alexander stratifications of character varieties, 
Ann. Inst. Fourier (Grenoble) {\bf 47}, (1997), 555--583.\\
{[23]} Y. Ihara, On Galois representations arising from towers of coverings of ${\bf P}^{1} \setminus \{ 0,1,\infty \}$, Invent. Math. {\bf 86}, (1986), 427-459.\\
{[24]} --------, Braids, Galois groups, and some arithmetic functions,
   In: Proc. of the International Congress of Mathematicians (Kyoto, 1990), 99-120.\\
{[25]} S. Iyanaga, T. Tamagawa, Sur la theorie du corps de classes sur le corps des nombres rationnels.  J. Math. Soc. Japan {\bf 3}, (1951), 220--227.\\
{[26]} K. Iwasawa, On ${\bf Z}_{l}$-extensions of algebraic number fields, Ann. of Math. (2) {\bf 98}, (1973), 246--326.\\
{[27]} T. Kadokami, Y. Mizusawa, Iwasawa type formulas for covers of a link in a rational homology sphere, J. Knot Theory and Its Ramification, {\bf 17}, No. 10 (2008), 1--23.\\
{[28]} M. Kapranov,  Analogies between the Langlands correspondence and topological quantum field theory, Progress in Math., {\bf 131}, Birkh\"{a}user, (1995), 119--151.\\
{[29]}--------, Analogies between number fields and 3-manifolds, unpublished note, Max-Planck, (1996).\\
{[30]} K. Kato, Lectures on the approach to Iwasawa theory for Hasse-Weil $L$-functions via $B_{dR}$, in: Arithmetic algebraic geometry, Lecture Note in Math., {\bf 1553}, Springer, (1993), 50-163.\\
{[31]} P. Kirk, C. Livingston, Twisted Alexander invariants, Reidemeister torsion, and Casson-Gordon invariants. Topology {\bf 38}, (1999), no. 3, 635--661.\\
{[32]} F. Knudsen, D. Mumford, The projectivity of the moduli
   space of stable curves. I. Preliminaries on "det" and "Div". Math. Scand. {\bf 39}, (1976), no. 1, 19--55.\\
{[33]} H. Koch, Galoissche Theorie der $p$-Erweiterungen, Springer, Berlin-New York, VEB Deutscher Verlag der Wissenschaften, Berlin, 1970.\\
{[34]} T. Kohno, Conformal field theory and topology,  Transl. of Math. Monographs, {\bf 210}, Iwanami Series in Modern Math., AMS, 2002.  \\
{[35]} M. Kontsevich, Vassiliev's knot invariants, Advances in Soviet Mathematics, {\bf 16}, (1993), 137-150. \\
{[36]} M. Kurihara, Iwasawa theory and Fitting ideals, J. Reine Angew.
   Math. {\bf 561}, (2003), 39-86.\\
{[37]} J. Labute, Mild pro-$p$ groups and Galois groups of $p$-extensions of ${\bf Q}$,  J. Reine Angew. Math. 596 (2006), 155--182. \\
{[38]} Tkhang Le, Varieties of representations and their subvarieties
   of homology jumps for certain knot groups, Russian Math. Surveys, {\bf 46},
(1991), 250--251.\\
{[39]}  W. L\"{u}ck, Analytic and topological torsion for manifolds with boundary and symmetry, J. Diff. Geo., {\bf 37}, (1993), 263-322.\\
{[40]} --------, $L^{2}$-invariants: theory and applications to
   geometry and $K$-theory, Ergebnisse der Mathematik und ihrer Grenzgebiete. 3. Folge. A Series of Modern Surveys in Math., {\bf 44}, Springer-Verlag, Berlin, 2002. \\
{[41]} T. Maeda, Lower central series of link groups, Ph.D. dissertation, Univ. Toronto, 1983.\\
{[42]} Yu.I. Manin, M. Marcolli, Holography principle and arithmetic of algebraic curves, Adv. Theor. Math. Phys. {\bf 5}, (2001), no. 3, 617--650.\\
{[43]} W. Massey, L. Traldi, On a conjecture of K. Murasugi, Pacific J. Math., {\bf 124}, (1986), 193-213.\\
{[44]} B. Mazur, Notes on \'{e}tale cohomology of number fields. Ann. Sci. \'{E}cole Norm. Sup. (4) {\bf  6}, (1973),  521-552. \\
{[45]} --------, Deforming Galois representations, In: Galois groups over ${\bf Q}$,  Math. Sci. Res. Inst. Publ., {\bf 16}, Springer, (1989), 385-437.\\ 
{[46]} --------, The theme of $p$-adic variation, Mathematics: frontiers and
perspectives, AMS, RI, (2000), 433--459.\\
{[47]} J. Milne, Arithmetic duality theorems, Perspectives in Math., {\bf 1}, Academic Press, Inc., Boston, MA, 1986. \\
{[48]} J. Milnor, Isotopy of links, In: Algebraic Geometry and Topology, A symposium in honour of S. Lefschetz (edited by R.H. Fox, D.S. Spencer and W. Tucker), Princeton Univ. Press, Princeton,  (1957), 280-306.\\
{[49]} --------, Infinite cyclic coverings, 1968 Conference on the Topology of Manifolds (Michigan State Univ., E. Lansing, Mich., 1967) pp. 115--133 Prindle, Weber \& Schmidt, Boston, Mass.\\
{[50]} M. Morishita, Milnor's link invariants attached to certain Galois groups over {\bf Q}, Proc. Japan Academy, {\bf 76}, No.2, (2000), 18-21.\\
{[51]} --------, On certain analogies between knots and primes, J. Reine Angew. Math. {\bf 550}, (2002) 141-167.\\
{[52]} --------, A theory of genera for cyclic coverings of links, Proc. Japan Academy, {\bf 77}, No.7, (2001), 115-118. \\
{[53]} --------, On capitulation problem for 3-manifolds, in: Galois Theory and Modular Forms, Dev. in Math., {\bf 11}, Kluwer, (2003), 305-313.\\
{[54]} --------, Milnor invariants and Massey products for prime numbers, Compositio Math., {\bf 140}, (2004), 69-83.\\
{[55]} --------, Analogies between knots and primes, 3-manifolds and number fields, preprint, (2003).\\
{[56]} --------, On analogies between primes and knots, (Japanese), Sugaku, {\bf 58}, No.1, (2006), 40--63.\\
{[57]} --------, On the Alexander stratification in the deformation space of Galois characters, Kyushu J. Math., {\bf 60}, (2006), 1-10.\\
{[58]} --------, Milnor invariants and $l$-class groups, in: Geometry and Dynamics of Groups and Spaces, In Memory of Alexander Reznikov, Progress in Math. Vol. {\bf 265}, Birk\"{a}user, 589-603, 2007.\\
{[59]} --------, Knots and Primes, (Japanese), Springer-Japan, 2009.\\
{[60]} M. Morishita, Y. Terashima, Arithmetic topology after Hida theory, in: Intelligence of Low Dimensional Topology 2006, World Scientific Publishing Co. in the Knots and Everything Book Series, vol {\bf 40}, 213-222.\\
{[61]} --------, Geometry of polysymbols, Mathematical Research Letters, {\bf 15}, 95-116, 2008.\\ 
{[62]} --------, Chern-Simons variation and Deligne cohomology, to appear in: Spectral Analysis in Geometry and Number Theory on Professor Toshikazu Sunada's 60th birthday,  Contemporary Math., AMS, 2008.\\
{[63]} K. Murasugi, Nilpotent coverings of links and Milnor's invariant, Low-dimensional topology, London Math. Soc. Lecture Note Ser., {\bf 95}, Cambridge Univ. Press, Cambridge-New York  (1985), 106-142.\\
{[64]} W. Neumann, D. Zagier, Volumes of hyperbolic three-manifolds, 
Topology {\bf 24}, (1985), no. 3, 307--332. \\
{[65]} A. Noguchi, A functional equation for the Lefschetz zeta functions of infinite cyclic coverings with an application to knot theory, Spring Topology and Dynamical Systems Conference. Topology Proc. {\bf 29} (2005), no. 1, 277--291.\\
{[66]} --------, Zeros of the Alexander polynomial of a knot, Osaka J. Math. {\bf 44}, (2007), no. 3, 567--577.\\
{[67]} R. Porter, Milnor's $\bar{\mu}$-invariants and Massey products, Trans. Amer. Math. Soc. {\bf 275}, (1980), 39-71.\\
{[68]} N. Ramachandran, A note on arithmetic topology,  C. R. Math. Acad. Sci. Soc. R. Can. {\bf 23}, (2001), no.4, 130--135. \\
{[69]} L. R\'{e}dei, Ein neues zahlentheoretisches Symbol mit Anwendungen auf die Theorie der quadratischen Zahlk\"{o}rper, I
 J. Reine Angew. Math. {\bf 180}, (1938), 1-43.\\
{[70]} A. Reznikov, Three-manifolds class field theory (Homology of coverings for a nonvirtually $b_{1}$-positive manifold), Sel. math. New ser. {\bf 3},  (1997), 361-399.\\
{[71]} ---------, Embedded  incompressible surfaces  and homology of ramified coverings of three-manifolds, Sel. math. New ser. {\bf 6}, (2000), 1--39.\\
{[72]} A. Reznikov, P. Moree, Three-manifold subgroup growth, homology of coverings and simplicial volume,  Asian J. Math. 1 (1997), no. 4, 764--768.\\
{[73]} M. Sakuma, Homology of abelian coverings of links and spatial graphs, 
Canad. J. Math. {\bf 47}, (1995), 201-224.\\
{[74}] A. Schmidt, Circular sets of prime numbers and $p$-extensions of the rationals, J. Reine Angew. Math. {\bf 596}, (2006), 115--130.\\
{[75]} P. Shalen, P.  Wagreich, Growth rates, $Z\sb p$-homology, and volumes of hyperbolic $3$-manifolds, Trans. Amer. Math. Soc. {\bf 331}, (1992), 895--917.\\{[76]} A. Sikora, Analogies between group actions on 3-manifolds and number fields, Comment. Math. Helv. {\bf 78}, (2003), 832--844.\\
{[77]} D. Silver, S. Williams, Mahler measure, links and homology growth. Topology {\bf 41}, (2002), 979--991.\\
{[78]} K. Sugiyama, An analog of the Iwasawa conjecture for a compact hyperbolic threefold, J. Reine Angew. Math. {\bf 613}, (2007), 35--50.  \\
{[79]} ---------, An analog of the Iwasawa conjecture for a complete hyperbolic threefold of finite volume, preprint.  \\
{[80]} T. Sunada, Fundamental groups and Laplacians, (Japanese), Kinokuniya, 1988.\\
{[81]} H. Suzuki, On the capitulation problem, In: Class field theory---its centenary and prospect (Tokyo, 1998), 483--507, Adv. Stud. Pure Math., {\bf 30}, 
Math. Soc. Japan, Tokyo, 2001. \\
{[82]} W. Thurston, The geometry and topology of 3-manifolds, Lect. Note, Princeton, 1977. \\
{[83]} L. Traldi, Milnor's invariants and the completions of link modules, 
 Trans. Amer. Math. Soc. {\bf 284}, (1984), 401-424.\\
{[84]} V.G. Turaev, Milnor's invariants and  Massey products, English transl. J. Soviet Math. {\bf 12}, (1979), 128-137 .\\
{[85]}--------, Nilpotent homotopy types of closed manifolds, Lecture Notes in Math., {\bf 1060}, Springer, (1984), 355-366.\\
{[86]} D. Vogel, On the Galois group of 2-extensions with restricted ramification, J. Reine Angew. Math. {\bf 581}, (2005), 117-150.\\
{[87]} J.-L. Waldspurger, Entrelacements sur Spec({\bf Z}), Bull. Sc. math. {\bf 100}, (1976), 113-139.\\
{[88]} Y. Yamamoto, Class number problems for quadratic fields (concentrating on the $2$-part). (Japanese) Sugaku, {\bf 40}, (1988), no. 2, 167--174.\\
\end{flushleft}
Graduate School of Mathematics, Kyushu University, Hakozaki, Fukuoka, 812-8581, Japan. morisita$@$math.kyushu-u.ac.jp}\\

\end{document}